\documentstyle[12pt,twoside,diagram]{article}
\font\twelvegtc=eufm10 scaled 1200

\font\ninegtc=eufm9
\font\sevengtc=eufm7

\newfam\gtcfam
\def\gtc{\fam\gtcfam\twelvegtc}
\textfont\gtcfam=\twelvegtc
\scriptfont\gtcfam=\ninegtc
\scriptscriptfont\gtcfam=\sevengtc
\font\twelveBBB=msbm10 scaled 1200
\font\tenBBB=msbm10
\font\sevenBBB=msbm7
\newfam\BBBfam
\def\BBB{\fam\BBBfam\twelveBBB}
\textfont\BBBfam=\twelveBBB
\scriptfont\BBBfam=\tenBBB
\scriptscriptfont\BBBfam=\sevenBBB
   \def\ZZ{{\BBB Z}}    \def\NN{{\BBB N}}
   \def\CC{{\BBB C}}    \def\PP{{\BBB P}}

\def\VEC#1,#2{(#1_1,#1_2,\dots,#1_{#2})}
\def\OVEC#1,#2{(#1_0,#1_1,#1_2,\dots,#1_{#2})}
\def\SET#1,#2{\{#1_1,#1_2,\dots,#1_{#2}\}}
\def\OSET#1,#2{\{#1_0,#1_1,\dots,#1_{#2}\}}
\def\FAM#1,#2{\ #1_1,#1_2,\dots,#1_{#2}\ }
\def\OFAM#1,#2{\ #1_0,#1_1,\dots,#1_{#2}\ }
\def\BSER#1,#2,#3{#1_0+#1_1#2+\cdots+#1_{#3}#2^{#3}}
\def\SER#1,#2{#1_0+#1_1#2+#1_2#2^2+\cdots}
\def\POL#1,#2,#3{#1_0#2^{#3}+#1_1#2^{#3-1}+\cdots+#1_{#3-1}#2+#1_{#3}}
\def\UPOL#1,#2,#3{#2^{#3}+#1_1#2^{#3-1}+\cdots+#1_{#3-1}#2+#1_{#3}}
\def\Ddots{\mathinner{\mkern1mu\raise1pt\hbox{.}\mkern2mu\raise5pt\hbox
                  {.}\mkern2mu\raise8pt\vbox{\kern7pt\hbox{.}}\mkern1mu}}
\def\bydef{\stackrel{\rm def}{=}}
\def\comp{{\scriptstyle\cdot}}
\def\bdot{{\scriptscriptstyle\bullet}}
\def\^#1#2{{^{#1}\!{#2}}}
\def\im{{\rm im}\,}
\def\coker{{\rm coker}\,}
\def\rk{{\rm rk}\,}

\def\Id{{\rm Id}}

\def\hom{{\rm Hom}}
\def\Hom{{\cal H\it om}}

\def\ext{{\rm Ext}}
\def\Ext{{\cal E\it xt}}
\def\End{{\rm End}}
\def\SL#1,#2{{\rm SL}_{#1}({#2})}
\def\GL#1,#2{{\rm GL}_{#1}({#2})}
\def\bl{\!\in\!}

\def\fain#1,#2{\;\forall\,{#1}\bl{#2}\;}
\def\map#1,#2{#1\longrightarrow#2}
\def\MAP#1,#2,#3{#1\,\colon\;\map#2,#3}
\def\TO#1{\setbox0=\hbox{\kern5pt$\scriptstyle{#1}$\kern9pt}
            \setbox1=\hbox to\wd0{\rightarrowfill}
            \;\vbox{\copy0\penalty10000\nointerlineskip\copy1}\;}
\def\FROM#1{\setbox0=\hbox{\kern9pt$\scriptstyle{#1}$\kern5pt}
            \setbox1=\hbox to\wd0{\leftarrowfill}
            \;\vbox{\copy0\penalty10000\nointerlineskip\copy1}\;}
\def\sp#1,#2{\langle#1,#2\rangle}
\def\cc#1,#2{\left[\frac{\textstyle #1}{\textstyle #2}\right]}
\newcounter{No}
\newcounter{SubNo}[No]
\newcounter{SubSubNo}[SubNo]
\renewcommand{\theNo}{\S\arabic{No}}

\newcounter{cond}

\def\No#1{\setcounter{equation}{0}\refstepcounter{No}
          \par\vspace{1cm}\par\noindent
          {\large\bf\theNo.\hspace{2pt}#1.}\par\nopagebreak[4]}
\def\SubNo#1{\refstepcounter{SubNo}
             \vspace{2ex}\par\noindent
            {\bf\arabic{No}.\arabic{SubNo}.\hspace{2pt}#1.}\hspace{1ex}}
\def\SubSubNo#1{\refstepcounter{SubSubNo}\vspace{1ex}\par\noindent
       {\bf\arabic{No}.\arabic{SubNo}.\arabic{SubSubNo}.}\hspace{2pt}{#1}}
\def\EF{\endgroup\par\vspace{1ex}\par}
\def\EP{\par\nopagebreak[3]\noindent$\Box$\par}

\def\Lm{\SubSubNo{LEMMA.\hspace{1ex}}\begingroup\sl}
\def\Th{\SubSubNo{THEOREM.\hspace{1ex}}\begingroup\sl}
\def\Cl{\SubSubNo{COROLLARY.\hspace{1ex}}\begingroup\sl}
\def\proof{\par\noindent{\sc Proof.}\hspace{1ex}}
\def\Cs{\setcounter{cond}{0}\begin{description}}

\def\EC{\end{description}}
\def\DBAS{\{\^2E,\^1E,\^0E\}}
\def\BAS{\{E_0,E_1,E_2\}}
\def\ct{c_{\rm top}}
\def\M{{\gtc M}}
\def\N{{\gtc N}}
\def\GR{{\gtc G}}
\def\Z{{\gtc Z}}
\def\V{{\BBB V}}
\def\O{{\cal O}}
\def\SS{{\cal S}}
\def\F{{\cal F}}
\def\C{{\cal C}}
\def\G{{\cal G}}
\def\U{{\cal U}}
\def\T{{\cal T}}
\def\X{{\cal X}}
\def\E{{\cal E}}
\def\nb#1#2{{\cal N}_{#1/#2}}
\def\SL#1,#2{{\sl SL}_{#1}({#2})}
\def\GL#1,#2{{\sl GL}_{#1}({#2})}
\def\PGL#1,#2{{\sl PGL}_{#1}({#2})}
\def\^#1#2{{^\times\!{#2}}_{#1}}
\let\s=\sigma
\let\g=\gamma
\let\L=\Lambda
\let\l=\lambda
\let\f=\varphi
\let\t=\tau
\begin{document}
\begin{center}
  \small\sc
  Independent University of Moscow\\
  Max-Planck Institute f\"ur Mathematik\\
  \vspace{.5cm}
  \LARGE\bf
   On top Chern classes of universal bundles \\ on moduli spaces
 of rank two coherent sheaves on the projective plane, or \\
 How to compute the correlation function \\
 in SYM $N=2$ $N_f = 4$ quantum field theory on
 $\CC \PP^2$.%
  \\ \vspace{.5cm}
  \Large\rm
  Alexei L. Gorodentsev\footnote{e-mail: \verb|gorod@itep.ru|}\\
  Maxim I. Leyenson\footnote{e-mail: \verb|leyenson@gmail.com|}\\
  \vspace{.5cm}
  \small\rm (12 April, 1995)
\end{center}
\vspace{1ex}\par
%
%
\No{Introduction}

\SubNo{What is this paper about}\label{start}
  Let $\M\bydef\M_{\PP_2}(2;-1;k)$ be the moduli space of $\mu$-stable
  torsion free coherent sheaves $F$ with $\rk(F)=2$, $c_1(F)=-1$, 
  $c_2(F)=k$ on $\PP_2=\PP_2(\CC)$. $\M$ is a smooth
  $4(k-1)$-dimensional projective variety and there exists {\it the
  universal\/} rank $2$ torsion free sheaf $\F$ on $\PP_2 \times \M$ 
  defined uniquely up to the twisting by the pull-back of an invertible 
  sheaf on $\M$ (see, for example, \cite[ch.II, \S4]{OSS}).

  We fix the standard notations for the projections
  \begin{equation}\label{pr}
  \begin{diagram}
  \node[2]{\PP_2\times\M}
          \arrow{se,l}{\pi}
          \arrow{sw,l}{p}
  \\
  \node{\PP_2}\node[2]{\M}
  \end{diagram}
  \end{equation}
  and put $\G\bydef R^1\pi_*\F$. The sheaf $\G$ is locally free and
  $\rk\G=(k-1)$. The fibers of $\G$ at the closed points $F\bl\M$ are
  isomorphic to $H^1(\PP_2;F)$. We call $\G$ {\it the universal bundle
  on\/} $\M$. 

  Let $\G^{\oplus4}\bydef\G\oplus\G\oplus\G\oplus\G$.
  The sequence of the topological constants
  $a_k\bydef c_{\rm top}(\G^{\oplus4})=c_{4k-4}(\G^{\oplus4})$, where 
  $k\bl\NN$, appears as the correlation function in the $N=2$, $N_f=4$
  supersymmetric Yang-Mills theories with $SO(3)$ gauge group.
 
  Namely, consider $X=\CC\PP_2$ as real smooth 4-manifold equipped
  with some Riemann metric $g$ and the standart ${\rm 
  Spin}^\CC$-structure $c = 3h$. Also let $F\TO{}X$ be a smooth 
  complex 2-dimensional vector bundle with $c_1(F)=-1$,  $c_2(F)=k$ 
  equipped with a $SO(3)$-connection $a$ and the corresponding twisted 
  Dirac operator 
  $\MAP D_a,\C^\infty(F\otimes W^+),{\C^\infty(F\otimes W^-)}$,
  where $W^\pm$ are the spinor bundles presented by ${\rm 
  Spin}^\CC$-structure $c$.

  It is easy to check that the correlation functions in the QFT with 
  four flawers, i.~e.  with four spinor fields 
  $$\MAP\phi_i,X,{F\otimes W^+}\;,\qquad i=1,2,3,4\;, 
  $$ 
  can be calculated correctly on the classical level, because the 
  $\beta$-function of this theory vanishes. Absolute minima of the  
  Yang-Mills Lagrangian are achieved exactly on the $(a,\psi_i)$ 
  that satisfy the equations 
  \begin{equation}\label{difur} 
  \left\{
  \begin{array}{rcl}
        F_a^+&=&\sum\limits_{i=1}^4(\phi_i\otimes\overline{\phi}_i)_{00}\\
        D_a\phi_i&=&0\;,\quad i=1,2,3,4
  \end{array} 
  \right.\quad,
  \end{equation}
  where $\MAP F^+_a,H^0(F),{H^0(\Omega^{2+}\otimes F)}$ is the 
  Hodge-selfdual part of the curvature of $a$. The $S^1$-action on 
  each of the coupled spinor fields can be killed by the same 
  normalisation of the detrminantal connections as in 
  \cite[${\rm n}^{\rm o}$ 2.2]{PiTy}. So, for general metric $g$ 
  on $X$ the equations (\ref{difur}) have only 
  a finite (up to $SO(3)$ gauge) number of solutions.  This number 
  $a_k$ (where $k=c_2(E)$) does not depend on a choice 
  of metric (if the choice is general enough). 

  To apply the geometry for the calculation of constants $a_k$, we have
  to use the Fubini-Study metric $g$. 
  It is not general: all spinor fields $\phi_i$, which satisfy 
  (\ref{difur}) must be zero in this case, because of positive scalar 
  curvature of $g$.  Hence, the equations (\ref{difur}) take the form 
  $F_a^+=0$, i.~e. define the usual instantons (antiselfdual 
  connections on $F$). By Donaldson's theorem, the instantons 
  (considered up to $SO(3)$ gauge) are in 1-1 correspondense with 
  holomorphic $\mu$-stable
  structures on $F$, and hence, they are parametrized by the
  Zarisski open subset of $\M(2,-1,k)$, which consists of all 
  locally free stable sheaves. So, in order to calculate the
  actual values of $a_k$ in terms of $\M$, we have to consider an
  {\it obstruction bundle\/}. Since the obstruction spaces for 
  the existence of non-zero spinor $\phi\in\ker D_a$ coinsides with 
  $\coker D_a=H^1(F)$, 
  the obstruction bundle is exactly the universal bundle $\G$ on 
  $\M(2,-1,k)$ and the number of solutions of (\ref{difur}) for 
  general metric equals $c_{\rm top}(\G^{\oplus4})$ (i.~e. the 
  number of points where four general obstructions vanish 
  simultaniously). 

  Constants $a_k$ play an important role in physic and the sum 
  $f(q)=\sum a_kq^k$ is waited by physicists to be equal 
  to the $q$-decomposition of a modular form, because of the
  physical S-duality conjecture (see \cite{VW,Wi}). To check such kind of 
  statements mathematically, we have to find a clear and not too hard way 
  for the calculation of the constants $a_k$.  Such a way is presented 
  in this paper.

\SubNo{Approach and results}\label{result}
  In \ref{redtograss} we give a direct geometrical construction 
  of the moduli space
  $\M_{\PP_2}(2;-1;k)$. In particular, this reduces the calculation of 
  $a_k$ to the problem that can be solved principally by the Schubert 
  calculus over a ring, which is given by some explicit generators and 
  relations.  For this aim we use the approach of G.Ellingsrud and 
  S.A.Str{\o}mme from \cite{ES2}. We extend this approach to the more 
  general {\it helices\/} on $\PP_2$ and adapt it for the sheaf $\G$.

  Namely, consider the {\it Kronecker moduli space\/} $\N=\N(3;k,k-1)$
  defined as the space of the stable orbits (i.e. the geometric factor)
  of the natural representation of the reductive group
  $$\Bigl(\GL{k},\CC\times\GL k-1,\CC\Bigr)/\CC\!\cdot\!\{\Id\times\Id\}
  $$
  in the vector space ${\CC^*}^{k}\otimes\CC^3\otimes\CC^{(k-1)}$. The
  variety $\N$ was studied in \cite{Dr} and its Chow ring has been
  calculated via generators and relations in \cite{ES3}.
  There exists the {\it universal Kronecker module\/}
  $\U_1^*\otimes\CC^3\otimes\U_2 $ over $\N$, where $\U_1$ and $\U_2$ are
  the natural universal vector bundles on $\N$ with $\rk\U_1=k$ and  
  $\rk\U_2=(k-1)$ (their Chern classes are the ring generators for
  $A^\bdot(\N)$, see  \cite{ES3}).

  Denote by
  ${\GR}={\sl Gr}\Bigl(\,k-1\,,\;\CC^3\otimes\U_1\,\Bigr)\TO{\rm pr}\N$
  the relative Grassmann parametrizing all $(k-1)$-subbundles
  in the bundle $\CC^3\otimes\U_1$ on $\N$. Let $\SS$ be the universal
  subbundle on ${\GR}$. The next theorem follows from what we prove in 
  \ref{redtograss}.

\markboth{{\sl Alexei L. Gorodentsev, Maxim I. Leyenson}}
{{\it How to calculate $N=2$, $N_f=4$ correlation function on $\PP_2$}}

\Th\label{th}
  The moduli space $\M=\M(2;-1;k)$ is isomorphic to the subvariety
  ${\Z}\!\subset\!{\GR}$ defined as the zero scheme of a section of the
  bundle $\CC^3\otimes\SS^*\otimes{\rm pr}^*\U_2$ on ${\GR}$. Under this
  isomorphism the universal bundle $\G$ on $\M$ is identified with the
  restriction of the universal bundle $\U_2$ onto ${\Z}$.
\EF
 
  In \S3 we study an action of the maximal torus $T\subset\PGL3,\CC$ 
  on $\M$. This action comes from a toric variety structure on $\PP_2$
  and the fixed points $E\in\M^T$ of this action are represented by the
  {\it toric sheaves\/} on $\PP_2$. Using general technique of 
  A.A.K1yachko 
  (\cite{K1,K2}), we enumerate all connected components $Y\subset\M^T$ 
  in some combinatorial terms.
 
  Namely (see details in \S3), each connected component $Y\subset\M^T$
  is given by

  \smallskip\begingroup\leftskip=1cm

  \noindent\llap{--}
  an ordered collection of three non-negative integers $d_0,d_1,d_2$ 
  such that the sum $d=d_0+d_1+d_2$ is bounded by $0\le d<c_2(E)$,

  \smallskip

  \noindent\llap{--}
  an ordered collection of three positive integers $a_0,a_1,a_2$, which 
  satisfy three triangle inequalities $a_i<a_j+a_k$ and the equality 
  $(-a_0+a_1+a_2)^2-4a_1a_2=1-4(c_2(E)-d)$,

  \smallskip

  \noindent\llap{--}
  an ordered collection of three pictures like shown in (\ref{picsum}) on 
  the page \pageref{picsum} (i-th picture, $i=0,1,2$, is obtained from 
  a pair of Young diagrams $(\l_i,\mu_i)$ such that $|\l_i|+|\mu_i|=d_i$;
  the diagrams are shifted with respect to each other by $a_j$ cells 
  in the horizontal direction and by $a_k$ cells in the 
  vertical direction).

  \smallskip\endgroup

  \noindent
  Two such data leads to the same $Y$ iff they have the same triples 
  of numbers $d_i$, $a_i$ and the same (up to a parallel translation) 
  triples of pictures. The figure contained between two polygonal 
  boundaries of the Young diagrams $\l_i$, $\mu_i$ in (\ref{picsum}) 
  splits into connected components. Let $N_i$ be the number of the 
  components, which are strictly contained inside the intersection 
  of two dotted right angles and let $N=N_0+N_1+N_2$. Then the 
  number of Young diagram pair's triples $\{(\l_i,\mu_i)\}_{i=0,1,2}$ 
  leading to the same triple of pictures is $2^N$ and the corresponding 
  component $Y\subset\M^T$, is isomorphic to 
  $$\underbrace{\PP_1\times\PP_1\times\cdots\times\PP_1}_N\;,
  $$ 
  where the multipliers $\PP_1$ are in the natural 1-1 correspondence with 
  the connected pices between the Yong diagram boundaries described above.

  Using geometrical construction of $\M$ presented in \S2 we fix some
  toric structure on the universal bundle $\G$ over $\M$ and get the
  corresponding families of toric structures on $E$ runing 
  through any connected 
  component $Y\subset\M^T$. This let us describe the toric character 
  decomposition of all restrictions $\G|_Y$ (see \ref{hdecomp}). 
  
  In \S4 we explain how to apply the Bott residue formula for the 
  calculation of the consitsnts $a_k$. Unfortunately, we do not get any
  general closed answer, which express all $a_k$ directly in terms 
  of $k$ and/or some recursive rules, because the combinatorial data 
  required to evaluate the denominator in the Bott
  formula is too complicated. 
  But we present an explicit description of this data sufficient for 
  the numerical calculation $a_k$ by a computer. 
  The first values of $a_k$ (with respect to the normalizations of 
  $\G$ and $\M$ given in \ref{th}) are the following\footnote{
  they are calculated using MAPLE V.3}:  
  $$\mbox{\begin{tabular}{|c||c|c|c|c|c|c|c|}
              \hline
              $k$&1&2&3&4&5&6&7\\
              \hline 
              $a_k$&0&0&0&13&729&85026&15650066\\
              \hline 
          \end{tabular} 
         }
  $$

\SubNo{About the helices on $\PP_2$}\label{abouthel}
   All what we need about the helices on $\PP_2$ can be read in
   \cite{GoPn,GoRu}. We will use the terminology and the notations from
   these papers. Everybody who is interesting only in the above
   numbers can suppose without any loss that the {\it helix foundations\/}
   $\BAS$ and $\DBAS$, which will be used in what follows, are equal to
   the triples
   $$\{\O(-1),\,\Omega_{\PP_2}(1),\,\O\}\qquad{\rm and}\qquad
   \{\O,\,\O(1),\,\O(2)\}
   $$
   of sheaves on $\PP_2$.

   The more detailed review of the helix theory and the list of references
   can be founded in \cite{GoRew}.

\SubNo{About toric varieties and toric bundles}
  All toric things what we need (with a lot of further references) can be founded
  in \cite{K1,K2}. Our terminology and notations will be very closed to 
  the ones used in these papers. For convenience of readers we recall some 
  basic facts adapted for our framework at the beginning of \S3.

\SubNo{Acknowledgements} We 
  are grateful to A.~N.~Tyurin who draws our attention to this subject and 
  explanes a lot of physical things, and also to 
  S.~A.~Kuleshov for very useful remarks and improvements.

  This paper was finished when the first of authors was staying at Max-Planck
  Institute f\"ur Mathematik in Bonn and he would like to express here
  his gratitude for the hospitality. 

\No{Geometrical description of moduli space}
\label{redtograss}

\SubNo{The bases used for the representation of sheaves on $\PP_2$}
  Let $\BAS$ be a helix foundation on $\PP_2$ such that $E_0=\O(-1)$ and
  $\mu(E_2)\ge\mu(F)=-1/2$. It follows from \cite{GoRu} that the pair
  $\{E_1,E_2\}$ in such a foundation is the pair of consequent elements
  $\{X_\nu,X_{\nu+1}\}$ in the infinite sequence of sheaves
  $\{X_\nu\}_{1\le\nu<\infty}$ defined recursively by the relations
  $X_1=\Omega_{\PP_2}(1)$, $X_2=\O$, and
  $$X_{\nu+2}=\coker\Biggl(X_{\nu}\TO {\rm coev}
                 \hom(X_{\nu},X_{\nu+1})^*\otimes X_{\nu+1}\Biggr)\;.
  $$
  There exists also the natural identification
  $$\hom(X_\nu,X_{\nu+1})=\left\{
    \begin{array}{cl}
       \V^*&\mbox{for even $\nu$}\\
       \V&\mbox{for odd $\nu$}
    \end{array}
    \right.\;,
  $$
  where $\V^*=H^0(\O_{\PP_2}(1))$. We denote the space $\hom(E_i,E_j)$ by
  $V_{ij}$. As we have just explained, $\dim V_{12}=3$.

  Let $\DBAS$ be the {\it left dual foundation\/}. It consists of the
  sheaves $\^2E=E_2$, $\^0E=E_0(3)=\O(2)$, and
  $$\^1E=R_{E_2}(E_1)={\rm coker}\Biggl(
                E_1\TO{\rm ev}\hom(E_1,E_2)^*\otimes E_2\Biggr)
  $$
  The left dual foundation is uniquely defined by the orthogonality
  conditions
  $$\ext^{(2-i)}(\^jE,E_k)=\left\{
    \begin{array}{cl}
          \CC&{\rm for}\;i=j=k\\
          0&\mbox{in the other cases}
    \end{array} \right.\;.
  $$
  If the pair $\{E_1,E_2\}$ coincides with the pair $\{X_\nu,X_{\nu+1}\}$ 
  of consequent elements from the sequence $\{X_\nu\}$ described above, then
  the dual pair $\{\^2E,\^1E\}$ coincides with the pair $\{X_{\nu+1},X_{\nu+2}\}$
  from the same sequence. So,  in terms of $X_\nu$, every pair of dual foundations, 
  with is appropriate for us, can be written as
  $$\{E_0,X_{\nu},X_{\nu+1}\}\qquad{\rm and}\qquad
    \{X_{\nu+1},X_{\nu+2},\^0E\}
  $$ 
  and there exists the natural coincidence
  $\hom(X_\nu,X_{\nu+1})=\hom(X_{\nu+1},X_{\nu+2})^*$.

\SubNo{Monad on $\PP_2$}
  For $i=0,1,2$ the inequalities $\mu(F)\le\mu(\^iE)\le\mu(F(3))$ between
  the Mumford slopes imply (via the stability and the Serre duality) that
  $$\ext^0(\^iE,F)=\ext^2(\^iE,F)=0\;.
  $$
  So, only $\ext^1(\^iE,F)$ can be non-zero. We denote this vector space
  by $U_i$ and put $u_i=\dim U_i$.
  The Beilinson spectral sequence (see \cite{GoPn}) shows that $F$
  is the cohomology sheaf of the monad
  \begin{equation}\label{m}
  0\TO{}U_0\otimes E_0\TO{\iota_1}U_1\otimes E_1\TO{\iota_2} U_2\otimes
  E_2\TO{}0\;.
  \end{equation}
  Note that $U_0=H^1(F(-2))$ and $u_0=k-1$. In general, it follows from
  Riemann-Roch that
  $$u_i=-\chi(\^iE,F)=\rk(\^iE)\cdot k+2c_1(\^iE)-m(\^iE)\;,
  $$
  where $m(\^iE)\bydef(1+c_1(\^iE)^2)/\rk(\^iE)$ is an integer (see
  \cite{GoRu}), which can be called the {\it Markov characteristic\/} of
  the exceptional vector bundle $\^iE$ on $\PP_2$. In particular, for the
  simplest pair of the dual foundations
  $\{\O(-1),\,\Omega_{\PP_2}(1),\,\O\}$ and $\{\O,\,\O(1),\,\O(2)\}$
  we have the dimensions $u_1=k$ and $u_2=k-1$.

 \SubSubNo{\bf Some arithmetical conditions on $u_i$.}\label{arprop}
  The classes $\{e_0,e_1,e_2\}$ of the sheaves $\BAS$ form an integer
  semiorthonormal basis of the Mukai lattice $K_0(\PP_2)$.  In terms of
  this basis, the class $f$ of $F$ is decomposed as
  $$f=-u_0e_0+u_1e_1-u_2e_2\;,
  $$
  where $u_i=\chi(\^ie,f)$ and the left dual basis $\{\^ie\}_{i=2,1,0}$
  satisfy the relations $\chi(\^ie,e_j)=(-1)^i\delta_{ij}$.
  If we use this decomposition in order to calculate $\chi(e_0,f)$, then we
  get
  $$\chi(e_0,f)=-u_0+\chi(e_0,e_1)\cdot u_1+\chi(e_0,e_2)\cdot u_2\;.
  $$
  and hence,
  $$\chi(e_0,e_1)\cdot u_1+\chi(e_0,e_2)\cdot u_2=\chi(e_0,f)-\chi(\^0e,f)
                                                 =\chi(e_0,f)-\chi(f,e_0)
    =3\delta(e_0,f)\;,
  $$
  where the determinant $\delta(X,Y)\bydef \rk(X)c_1(Y)-\rk(Y)c_1(X)$ 
  coincides by Riemann-Roch with the skew-symmetric part of 
  $\frac{\textstyle 1}{\textstyle 3}\chi(X,Y)$.
     
  On the other side, standard arithmetical properties of exceptional pairs
  on $\PP_2$ (see \cite{GoRu}) imply the identities 
  $$\chi(e_0,e_1)=3\delta(e_0,e_1)=3\rk(E_2)\;.
  $$
  $$\chi(e_0,e_2)=3\delta(e_0,e_2)=3\rk(\^1E)\;,
  $$
  Hence, we have the identity
  \begin{equation}\label{nodrel}
  \rk(E_2)\cdot u_1-\rk(\^1E)\cdot u_2=-\rk(E_0)-2c_1(E_0)\;.
  \end{equation}
  Since the right side equals 1 for $E_0=\O(1)$, in this case 
  $u_0=\dim U_0$ and $u_1=\dim U_1$ are coprime.
  The other thing (which will be used below to study the stability) is that
  $$\delta(E_0,F)=\delta(E_0,K)=1\;,
  $$
  where $K=\ker(\iota_2)$ is the kernel of the
  second map from the monad (\ref{m}). 

\SubNo{Monad on $\PP_2\times\M$}
  The machinery for the Beilinson decomposition of coherent sheaves on
  $\PP_n$ in terms of dual helix foundations has the straightforward
  extension (see \cite{Or}) to the relative case, 
  where we consider the projectivization
  of a vector bundle over a base variety (instead of $\PP_n$
  over a field).  In particular, any coherent sheaf
  $\X$ on $\PP_2\times\M$ has a natural representation as the limit
  of the following Beilinson spectral sequence
  $$E_1^{\alpha,\beta}=\pi^*\Ext^\beta_\pi(p^*(\^{2+\alpha}E),\X)\otimes
    p^*E_{2+\alpha}\;;\quad\alpha=-2,-1,0\;,\quad\beta=0,1,2\quad,
  $$
  where
  $\Ext^\beta_\pi(p^*(\^{2+\alpha}E),\X)\bydef
  R^\beta\pi_*\Hom(p^*(\^{2+\alpha}E),\X)$ and $p$, $\pi$ are the
  projections from (\ref{pr}).

  So, the universal bundle $\F$ on $\PP_2\times\M$ is the cohomology sheaf
  of the monad
  \begin{equation}\label{M}
  0\TO{}\pi^*\F_0\otimes p^*E_0\TO{}
          \pi^*\F_1\otimes p^*E_1\TO{}
          \pi^*\F_2\otimes p^*E_2\TO{}0\;,
  \end{equation}
  where we denote by $\F_i$ the sheaves $\Ext^1_\pi(p^*(\^iE),\F)$ on
  $\M$. In the fibers at the points of $\M$ the monad (\ref{M}) induce
  various monads (\ref{m}) on $\PP_2$.

\SubNo{Universal complex on $\M$}
  The complex of the locally free sheaves (\ref{M}) is {\it adapted\/}
  (in the sense of \cite[Ch III$,\,{\rm n}^o\,6.2$]{GeMa}) to calculate the
  derived functor $R\pi_*R\Hom(p^*E_0, \F)$, because
  $$R^q\Hom(p^*E_0,\,\pi^*\F_i\otimes p^*E_i)=0\mbox{\ for $q>0$}\;;
  $$
  $$R^q\pi_*\Hom(p^*E_0,\,\pi^*\F_i\otimes p^*E_i)= \F_i^*\otimes
  R^q\pi_*\Hom(p^*E_0,\,p^*E_i)=
             \left\{\begin{array}{cl}
                 V_{0i}\otimes\F_i^*&\mbox{for $q=0$}\\
                 0&\mbox{for $q\ne0$}
             \end{array}\right. \;.
  $$
  Hence, applying the functor $R\pi_*R\Hom(\,p^*E_0,\,?\,)$ to
  the complex (\ref{M}), we get the following complex of locally free
  sheaves on $\M$:
  \begin{equation}\label{MM}
  0\TO{}\F_0\TO{}V_{01}\otimes\F_1\TO{}V_{02}\otimes\F_2\TO{}0\;.
  \end{equation}
  This triple is exact on the left term and has on the middle and on
  the right the cohomology sheaves $\Ext^0_\pi(p^*E_0,\F)$
  and $\Ext^1_\pi(p^*E_0,\F)$. Both cohomologies can be not
  locally free\footnote{Their fibers at the point $F\bl\M$ are
  $H^0(F(1))$ and $H^1(F(1))$ and these spaces can jump.}.

  Note that for the simplest triple 
  $\{E_0,E_1,E_2\}=\{\O(-1),\Omega(1),\O\}$
  the sheaf $\F_2$ in the right term of (\ref{MM}) is exactly the universal
  sheaf $\G$ on $\M$.

   We are going to
  construct the closed embedding of $\M$ into the relative Grassmann
  ${\GR}$ over the Kronecker moduli space. Under this embedding the
  universal complex (\ref{MM}) will coincide with the restriction of some
  natural universal sequence of bundles on the Grassmann ${\GR}$.

\SubNo{The map from $\M$ to $\N(V_{01};U_0,U_1)$}
  The surjection $\iota_2$ from the monad (\ref{m}) gives a {\it
  Kronecker module\/}, i.~e. a tensor $\kappa_F$ in the space
  \begin{equation}\label{space}
  U_1^*\otimes V_{12}\otimes U_2\;.
  \end{equation}

  In accordance with J.-M.Drezet's results (see \cite[Prop.27]{Dr}) 
  the tensor $\kappa_F$ is 
  stable with respect to the natural action of the group
  \begin{equation}\label{group}
  G={\rm GL}(U_1)\times{\rm GL}(U_2)/\CC\!\cdot\!\{\Id\times\Id\}\;,
  \end{equation}
  if and only if the corresponding map $\iota_2$ is surjective and its kernel $K$ 
  from the exact triple
  $$0\TO{}K\TO{}U_1\otimes E_1\TO{\iota_2}U_2\otimes E_2\TO{}0
  $$
  is stable sheaf on $\PP_2$. In order to check that these conditions hold
  in our case, we use the other exact triple induced by the monad (\ref{m}):
  $$0\TO{}U_0\otimes E_0\TO{}K\TO{}F\TO{}0\;.
  $$
  and note that $K$ belongs to the category spaned by $E_1$ and $E_2$,
  i.~e. $\hom(K,E_0)=0$. So, by \ref{arprop}
   we are in a position to apply the following general lemma.

\Lm\label{stabext}
  Let $E$ be a stable locally free sheaf  and $F$ -- any sheaf such that 
  $$\delta(E,F)=\rk(E)c_1(F)-\rk(F)c_1(E)=1\;.
  $$
  If an extension $K$ of the form
  $$0\TO{}V\otimes E\TO{}K\TO{}F\TO{}0\;
  $$
  satisfy the condition $\hom(K,E)=0$, then $F$ is Mumford-stable 
  if and only if $K$ is.
\EF

\proof
  In terms of an integer 2-dimensional lattice generated by additive functions
  $(\rk, c_1)$, the determinantal conditions $\delta(E,K)=\delta(E,F)=1$ 
  mean that there are no integer points except for vertices in two parallelograms 
  spaned by the classes of $E$ and $K$ and by the classes of $E$ and $F$.
  In particular, any class $X$, which satisfy the conditions
  $$\mbox{ either }\quad\left\{
    \begin{array}{ccc}
    \mu(X)&<&\mu(K)\\
    \rk(X)&<&\rk(K)
    \end{array}
    \right.\quad\mbox{ or }\quad\left\{
    \begin{array}{ccc}
    \mu(X)&<&\mu(F)\\
    \rk(X)&<&\rk(F)
    \end{array}
    \right.  
  $$
  have to satisfy either the inequality $\mu(X)<\mu(E)$ or the system of 
  equalities
  $$\left\{
    \begin{array}{ccc}
       \rk(X)&=&m\cdot\rk(E)\\
       c_1(X)&=&m\cdot c_1(E)
    \end{array}\right.
  $$ 
  for some integer $m$.

  Now we consider consequently both implications of the lemma.
  
  Let $K$ be stable. If $F$ is not and admits a torsion free factor
  $F\TO{}Q\TO{}0$ with $\mu(Q)\le\mu(F)$ and $\rk(Q)<\rk(F)$, then 
  automatically $\mu(Q)\le\mu(E)<\mu(K)$.
  This is impossible, because $Q$ is the factor-sheaf for $K$ too.

  Conversely, let $F$ be stable and $K$ be not and let $K\TO{}G\TO{}0$ be 
  the minimal stable torsion free  Harder-Narasimhan's factor with 
  $\mu(G)\le\mu(K)$ and $\rk(G)<\rk(K)$. If $\hom(E,G)=0$ then 
  $\hom(F,G)\ne0$ and  we get the contradiction as above. If there exists a 
  non-zero map $E\TO{\psi}G$, then $\mu(G)=\mu(E)$, because $V\otimes E$
  is semistable. Since $E$ is stable, locally free, and has the same slope
  as $G$, the
  non-zero map  $E\TO{\psi}G$ is an isomorphism. So, there exists a non-zero 
  map from $K$ to $E$.  This contradicts to the last assumption of the lemma.

\EP

\Cl \label{MinN}
  Taking $F\longmapsto\kappa_F$, we define a family of stable
  $G$-orbits of the Kronecker modules. Since the family is parametrized by
  the points of $\M$, this gives an algebraic morphism $\MAP\kappa,\M,\N$,
  where $\N\bydef\N(V_{12};U_1,U_2)$ is the space of stable orbits (i.e.
  the geometric factor) of the natural representation of the group $G$
  from (\ref{group}) in the vector space (\ref{space}).  
\EF

\EP

  Since the dimensions $u_1=\dim U_1$ and $u_2=\dim U_2$ are coprime,
  all semistable tensors in (\ref{space}) are stable
  (see \cite{Dr}), the action of $G$ on the set of stable tensors is free
  (see \cite{ES3}), and the spaces
  \begin{equation} \label{gmod}
    U_1\otimes\biggl(\det(U_1)^{-\nu_1}\otimes\det(U_2)^{-\nu_2}\biggr)^{\dim U_1}
    \qquad{\rm and}\qquad
    U_2\otimes\biggl(\det(U_1)^{-\nu_1}\otimes\det(U_2)^{-\nu_2}\biggr)^{\dim U_2}
  \end{equation}
  admit the natural structures of $G$-modules (if $\nu_1$ and $\nu_2$ are 
  such that $\nu_1u_1+\nu_2u_2=1$).
  Hence, $\N$ is a smooth projective variety and $G$-modules
  (\ref{gmod}) induce some vector bundles on $\N$. We denote these
  bundles by $\U_1$ and $\U_2$. Note that by the construction we have the 
  natural line bundle's isomorphism $\det\U_1 = \det\U_2$.

  The universal Kronecker module $\U_1^*\otimes V_{12}\otimes\U_2$ has the
  tautological global section, which can be considered as the homomorphism
  $\U_1\TO{}V_{12}\otimes\U_2$. This homomorphism induces the map
  $$V_{01}\otimes\U_1\TO{}V_{01}\otimes V_{12}\otimes\U_2\;,
  $$
  which can be composed with the multiplication map
  $$V_{01}\otimes V_{12}\otimes\U_2
              \TO{\mu\otimes\Id}V_{02}\otimes\U_2\;.
  $$
  Denote the resulting map by
  \begin{equation}\label{umk}
  V_{01}\otimes\U_0\TO{\gtc i}V_{02}\otimes\U_2
  \end{equation}
  The next corollary follows immediately from the construction.

\Cl
   The homomorphism $V_{01}\otimes\F_1\TO{}V_{02}\otimes\F_2$
   from (\ref{MM}) coincides with the inverse image of (\ref{umk})
   via the map $\MAP\kappa,\M,\N$ from \ref{MinN}

\EF
\EP

  The corollary explains why the map $\MAP\kappa,\M,\N$ is not an
  isomorphism: the kernel of the homomorphism
  $V_{01}\otimes\F_1\TO{}V_{02}\otimes\F_2$ from (\ref{MM}) must
  contain  the subbundle $\F_0\subset V_{12}\otimes\F_1$. To take this
  into account we have to lift (\ref{umk}) on the Grassmann
  parametrizing all the subbundles in question.

\SubNo{The embedding $\M\hookrightarrow{\GR}$}
  Let ${\GR}\bydef{\sl Gr}\Bigl(\,u_0\,,\;V_{01}\otimes\U_1\,\Bigr)\;
  \TO{\rm pr}\;\N$ be the relative Grassmann parametrizing the rank
  $u_0=\dim U_0=\rk\F_0$ subbundles in the vector bundle
  $V_{01}\otimes\U_1$ on $\N$. Denote by
  $$0\TO{}\SS\TO{\gtc j}V_{12}^*\otimes{\rm pr}^*\U_1
  $$
  the canonical inclusion of the universal subbundle $\SS$ on
  ${\GR}$. Let ${\Z}\subset {\GR}$ be the zero scheme of the composition
  \begin{equation}\label{zeromap}
  \SS\TO{\gtc j}V_{01}\otimes{\rm pr}^*\U_1\TO{{\rm pr}^*({\gtc i})}
  V_{02}\otimes{\rm pr}^*\U_2
  \end{equation}

\Lm\label{ident}
  The map $\MAP\kappa,\M,\N$ from \ref{MinN} lifts to a closed embedding
  $\MAP{\bf k},\M,{\GR}$, which maps $\M$ isomorphicaly onto ${\Z}$. Under
  this map the complex (\ref{MM}) on $\M$ is identified with the
  restriction of the sequence (\ref{zeromap}) onto ${\Z}$. 

  In particular, for the simplest triple
  $\{E_0,E_1,E_2\}=\{\O(-1),\Omega(1),\O\}$
  the restriction of the bundle ${\rm pr}^*\U_2$ onto ${\Z}$ 
  coincides with the universal bundle $\F_2=\G$ .
\EF

\proof
  It is sufficient to prove that the fiber of the triple 
  (\ref{zeromap}) over any closed point of $\Z$ induces on $\PP_2$ a monad 
  of the form (\ref{m}) with a $\mu$-stable cohomology sheaf.

  If we denote the fibers of the bundles $\SS$, ${\rm pr}^*\U_0$,
  and ${\rm pr}^*\U_1$ by $U_0$, $U_1$, and $U_2$, then for any
  closed point of $\Z$ the maps ${\gtc j}$ and ${\rm pr}^*({\gtc i})$
  give two Kronecker modules
  $$\kappa_1\in U_0^*\otimes V_{01}\otimes
    U_1\qquad{\rm and}\qquad \kappa_2\in U_1^*\otimes V_{12}\otimes U_2\;,
  $$
  which induce the sequence of sheaf homomrphisms on $\PP_2$:
  \begin{equation}\label{testmon}
  0\TO{}U_0\otimes E_0\TO{\iota_1}U_1\otimes E_1
  \TO{\iota_2}U_2\otimes E_2\TO{}0\;.
  \end{equation}
  Certainly, we have $\iota_2\comp\iota_1=0$, because this composition considered
  as a tensor from $V_{02}\otimes\hom(U_0,U_2)$ coincides with the composition
  $$U_0\TO{\kappa_1}
    U_1\otimes V_{01}\TO{\kappa_2\otimes\Id}
    U_2\otimes V_{12}\otimes V_{01}\TO{\Id\otimes\mu}U_2\;,
  $$
  which vanishes over $\Z$ by the construction.
  Further, the second Kronecker module $\kappa_2$ is stable. Hence, 
  by J.-M.Drezet \cite[Prop.27]{Dr} it induces the exact triple
  \begin{equation}\label{stabC}
  0\TO{}K\TO{}U_1\otimes E_1\TO{\iota_2}U_2\otimes E_2\TO{}0
  \end{equation}
  on $\PP_2$ and the kernel $K$ of this triple is stable. So, as soon 
  as we prove that $\iota_1$ is an monomorphism we get that the 
  complex (\ref{testmon}) is a monad and its homology sheaf 
  is stable by \ref{stabext}.   

  Since $K$ and $E_0$ are stable and $\delta(E_0,K)=1$, the image  
  $I=\im\biggl(\,U_0\otimes E_0\TO{\iota_1}K\,\biggr)$ 
  is semi-stable and has the same slope as $E_0$. Hence, all 
  Jordan-H\"older factors of $I$ are equal to $E_0$ and
  $I$ is a direct sum of some copies of $E_0$, because $E_0$ 
  is exceptional. In other words, $\ker(\iota_1)=W\otimes E_0$ 
  for some subspace $W\subset U_0$.  
  \begin{equation}\label{matraz}
  \begin{diagram}
  \node[2]{0}
  \node[2]{0}
  \\
  \node{0}
          \arrow{e}
  \node{\im(\iota_1)}
          \arrow{n}
          \arrow[2]{e}
  \node[2]{U_1\otimes E_1}
          \arrow{n}
  \\
  \node{0}
          \arrow{e}
  \node{U_0\otimes E_0}
          \arrow[2]{e,t}{\kappa_1\otimes\Id}
          \arrow{n}
          \arrow{ene,l}{\iota_1}
  \node[2]{\hom(E_0,U_1\otimes E_1)\otimes E_0}
          \arrow{n,r}{\rm ev}
  \\
  \node{0}
          \arrow{e}
  \node{W\otimes E_0}
          \arrow{n}
          \arrow[2]{e}
  \node[2]{U_1\otimes L_{E_0}E_1}
          \arrow{n}
  \\
  \node[2]{0}
          \arrow{n}
  \node[2]{0}
          \arrow{n}
  \end{diagram}
  \end{equation}
  Since the first Kronecker module $\kappa_1$ is injective, we can draw
  the commutative diagram shown on (\ref{matraz}).
  The right column of this diagram is the canonical exact triple
  \footnote{In the terminology of the helix theory, this triple
  define the left mutation of $E_1$ by $E_0$.}:
  $$0\TO{}E_1^\times\TO{}\hom(E_0,E_1)\otimes E_0\TO{\rm ev}E_1\TO{}0\;,
  $$
  tensored by the vector space $U_1$. The bottom row is 
  induced by the rest part of the diagram and shows that $W=0$,
  because $\hom(E_0,E_1^\times)=0$. So, $\ker(\iota_1)=0$ and the proof
  is finished.

\EP

\SubNo{How to calculate $a_k$ via Schubert}
  Extracting $\dim\M$ from the scalar square of any class $f\bl K_0(\PP_2)$
  coming from $F\bl\M$ (see \ref{arprop}),  we get 
  $$\dim\M=1-u_0^2-u_1^2-u_2^2+h_{01}u_0u_1+h_{12}u_1u_2-h_{02}u_0u_2\;,
  $$
  where $h_{ij}=\dim V_{ij}$. Since 
  $\dim\N=h_{01}u_0u_1-u_0^2-u_1^2+1=\dim\M+u_2(h_{02}u_0-h_{12}u_1+u_2)$,
  we have
  $$\dim\GR=u_2(h_{12}u_1-u_2)+\dim\N=\dim\M+h_{02}u_0u_2\;,
  $$
  i.~e. the codimension of $\M=\Z$ in $\GR$
  $$d\bydef h_{02}u_0u_2=3\rk(\^1E)=3(k-1)(u_1\rk(E_0)-1)\;.
  $$
  is equal to the rank of the bundle
  $\Hom(\,\SS\,,\;V_{02}\otimes{\rm pr}^*\U_2\,)$. Hence, $\Z$ localizes
  the top Chern class of this bundle (see \cite[Prop.14.1]{Fu}).
  Since the rational and numerical equivalences of cycles on $\M$ coincide
  to each other (see \cite{ES1}), we have $A^*(\M)={\gtc c}_d\cdot A^*(\GR)$,
  where ${\gtc c}_d\bl A^*(\GR)$ is the top Chern class of the bundle
  $\SS^*\otimes V_{02}\otimes{\rm pr}^*\U_2$. Hence, the constants $a_k$ can
  be calculated inside the ring $A^*(\GR)$. For the simplest helix foundation
  we get in this way the formula:
  $$a_k=c_{(k-1)^2}(\SS)^3\cdot c_{k-1}(\U_2)^4\;.
  $$
  The ring $A^*(\GR)$ is the {\it Schubert calculus\/} (see \cite[\S14.7]{Fu})
  over the ring $A^*(\N)$.  An explicit description of the last one via
  generators and relations can be found in \cite[Th.(6.9)]{ES3}.

\SubNo{Normalization of the universal bundle $\G$}
  Note that the identification of the bundles $\F_1$, $\F_2$ on $\M$ with 
  the bundles $\U_1$, $\U_2$ via the embeddind $\M\hookrightarrow\GR$ gives
  us a normalization of the universal bundle, i.~e. it fixes the 
  $\det\G\in{\sl Pic}(\M)$ uniquely. In fact, if we multiply the universal sheaf
  $\F$ on $\PP_2\times\M$ by $\pi^*{\cal L}$, where ${\cal L}\in{\sl Pic}(\M)$,
  then $\det\F_1$ and $\det\F_2$ are shifted in ${\sl Pic}(\M)$ by 
  ${\cal L}^{\otimes u_1}$ and by ${\cal L}^{\otimes u_2}$ respectively. 
  Since $u_1$ and $u_2$ are coprime, there is at most one shift, which 
  leads to the identity $\det\F_1=\det\F_2$. But we have 
  $\det\F_1=\det\U_1=\det\U_2=\det\F_2$ by the construction.

  If we take different basic helix foundations $\BAS=\{\O(-1),X_\nu,X_{\nu+1}\}$,
  then we get a sequence of different normalizations of the universal bundle
  on $\M$. Numerical experiments show that these normalizations leads to the
  different values for the constants $a_k$. In the calculations, which will be 
  present in the following two paragraphs, we will always suppose that the
  universal bundle $\G$ is normalized via the simplest foundation
  $\BAS=\{\O(-1),\Omega(1),\O\}$.


\No{Torus action on \M(-1,2,k)}

\SubNo{Toric structure on $\PP_2$}
  Let us fix some homogeneous coordinates $(x_0\colon x_1\colon x_2)$
  on $\PP_2$ and consider $\PP_2$ as a toric variety with respect to
  the action of the diagonal subtorus $T\subset{\rm PGL}_3(\CC)$ on this
  coordinates. We will represent elements of the torus by matrices of the
  form 
  $$\left(\begin{array}{ccc}
                        1&&\\&t_1^{m_1}\\&&t_2^{m_2}
          \end{array}\right)
  $$ 
  and often we will use $(m_1,m_2)$ as the Cartesian coordinates on the torus
  character's lattice $\Lambda(T)$.

 It is convenient to introduce the following three one-parameter 
 subgroups $\t_i\in\L(T)^*$:
  $$
  \t_0(t)=\left(
       \begin{array}{ccc}
             1&&\\&t^{-1}\\&&t^{-1}
       \end{array}\right)
       \,,\quad
  \t_1(t)=\left(
       \begin{array}{ccc}
             1&&\\&t\\&&1
       \end{array}\right)
       \,,\quad
  \t_2(t)=\left(
       \begin{array}{ccc}
             1&&\\&1\\&&t
       \end{array}\right)
       \quad,
  $$
  and use sometimes the trigonal coordinates $(m_0.m_1,m_2)$ on $\L(T)$,
  where
  $$m_i(\chi)\bydef\sp\t_i,\chi\;.
  $$
  Of course, here $(m_1,m_2)$ are the same as above and $m_0+m_1+m_2=0$. 

  Standard affine card $U_i=\{x_i\ne0\}\subset\PP_2$ coincides with
  ${\rm Spec}\CC[\s_i]$, where $\s_i\in\L(T)$ are some angles such that in
  Cartesian coordinates $(m_1,m_2)$ 
  \begin{eqnarray*}
    \s_0&\mbox{generated by the characters}& (1,0)\mbox{ and }(0,1)\\ 
    \s_1&\mbox{generated by the characters}& (-1,0)\mbox{ and }(-1,1)\\ 
    \s_2&\mbox{generated by the characters}& (-1,1)\mbox{ and }(0,-1)\;, 
  \end{eqnarray*}
  Characters, which are regular on 
  $U_{ij}\bydef U_i\cap U_j$, 
  form the following halfplanes $\s_{ij}\in\L(T)$:
  \begin{eqnarray*}
    \s_{12}&=& \{(d_1,d_2)\in\L(T)\,\colon\;m_0=-m_1-m_2\ge0\}\\ 
    \s_{02}&=& \{(d_1,d_2)\in\L(T)\,\colon\;m_1\ge0\}\\ 
    \s_{01}&=& \{(d_1,d_2)\in\L(T)\,\colon\;m_2\ge0\}\quad. 
  \end{eqnarray*}
  We will always identify the torus $T$ with the open dense orbit 
  $U_{012}\bydef U_0\cap U_1\cap U_2={\rm Spec}(\CC[\L(T)])$ by puting $t\bl T$
  into $tp\bl U_{012}$, where $p=(1\colon1\colon1)\in U_{012}$.

\SubNo{Torus action on sheaves}
  The torus $T$ acts on the set of classes of isomorphic torsion free sheaves
  on $\PP_2$ via
  \begin{equation}\label{tstar}
  \MAP t,F,{t_*(F)}\;.
  \end{equation}
  So, we have a torus action on $\M=\M(2,-1,k)$ and also on the set 
  of sheaves on $\M$. We are going to apply to this action the Bott 
  residue formula (see \cite{Bt}). In order to do this, we have to 
  describe the connected components of the fixed point locus 
  $\M^T\subset\M$ and to fix a {\it toric structure\/} on the universal 
  bundle $\G$.

\SubSubNo{\bf Toric structures.}
  A sheaf $E$ represent a point of $\M^T$ iff there exists a
  a collection of isomorphisms
  $$\MAP\f_t,t_*(E),E\qquad\mbox{for every }\; t\in T\;.
  $$
  Such a collection $\f$ is called a {\it toric structure\/} 
  on a sheaf $E$ if it satisfy the following additional condition:
  $\forall\,s,t\in T$ the natural diagram
  $$
  \begin {diagram}
    \node{(t\comp s)_*(E)}
        \arrow{s,l}{\Id}
        \arrow{e,t}{\f_{t\comp s}}
    \node{E}
    \\
    \node{t_*(s_*(E))}
        \arrow{e,b}{t_*(\f_{s})}
    \node{t_*(E)}
        \arrow{n,r}{\f_{t}}
  \end{diagram}
  $$
  is commutative. A sheaf equipped with a toric structure is called a 
  {\it toric sheaf\/}.

  The obstruction to the commutativity of the above diagram is
  represented by a cocycle on $T$ with values in $\hom(E,E)$ 
  (see \cite{BoKa}). In our case this cocycle vanishes, because $E$ is
  stable (and $\hom(E,E)=\CC$). So, every sheaf $E$, which represents a point
  of $\M^T$, admits a toric structure.
  
  Certainly, there are many different toric structures on the same 
  sheaf $E\bl\M$. For example, toric structures on the structure sheaf 
  $\O$ are parametrized by the torus characters: the structure 
  $\tilde{\chi}$, which corresponds to the character $\chi$, takes 
  $f\in\O(t^{-1}U)$ to $\tilde{\chi}_tf\in\O(U)$ defined by
  $$\hat{\chi}_tf(x)\bydef\chi(t)\cdot f(t^{-1}x)\;.
  $$
  Tensor multiplication of a toric sheaf $E$ by the structure sheaf $\O$ 
  equipped with a toric structure $\tilde{\chi}$ is called a 
  {\it shift of a toric structure on $E$ by the character $\chi$\/}.
  Any two different toric structures on an arbitrary torsion free toric sheaf 
  $E$ can be obtained from each other by the shift by a character.
 
\SubSubNo{\bf Uniformal fixation of the toric structures.}\label{tsf} 
  First of all, on the sheaves $\{\O,\O(1),\O(2)\}$ we fix the 
  tautological toric structures coming from the
  toric variety structure on $\PP_2=\PP\V$.  In these structures $T$-modules
  $H^0(\O)$, $H^0(\O(1))$, and $H^0(\O(2))$ coincide with the trivial
  1-dimensional $T$-module $\CC$ and with the standard representations
  of $T$ in $\V^*$ and $S^2\V^*$ respectively.

  Since the mutations of exceptional bundles are $T$-equivariant, the 
  previous toric structures induce a toric structure on each exceptional 
  vector bundle $E$ on $\PP_2$. In particular, we get the toric structures 
  on the basic sheaves $\{\O(-1),\Omega(1),\O\}$. The corresponding $T$-modules
  \begin{eqnarray*}
   V_{01}&=&\hom(\O(-1),\Omega(1))\\
   V_{12}&=&\hom(\Omega(1),\O)\\
   V_{02}&=&\hom(\O(-1),\O)
   \end{eqnarray*}
  coincide with standard representations of $T$ in the spaces $\V$, $\V$,
  and $\L^2(\V)$ respectively.
  
  So, we get a natural torus action on all geometrical objects used in \S2:
  $T$ acts on the Kronecker moduli space $\N(V_{12};U_1,U_2)$ (because 
  it acts on $V_{12}$), there is a natural toric structure on the bundles 
  $\U_1$, $\U_2$ (induced by the trivial torus action on $G$-modules 
  (\ref{gmod})), and finally, the torus acts on the Grassmanian $\GR$. 
  The last action preserves the subvariety 
  $\M=\M(2,-1,k)\subset\GR$ and the restriction of this action onto
  $\M$ coincides with the action (\ref{tstar}). Hence, we have a toric 
  structure on the universal bundle $\G=\U_2|_{\M}$ and this structure 
  leads to a specific toric structure on each toric bundle $E$ on $\PP_2$. 

  In order to describe the last one, remember that by (\ref{gmod}) there 
  is $T$-equivariant isomorphism $\det{\rm pr}^*\U_2=\det{\rm pr}^*\U_1$. 
  Taking its fiber over $E\bl\M^T\subset\GR$, we see that 
  $T$-modules $\det H^1(E)$ and $\det H^1(E(-1))$ must be isomorphic to each 
  other. In other words, the sum of all characters from $H^1(E)$ is equal 
  to the sum of all characters from $H^1(E(-1))$. It is known (see 
  \cite{K1,K2}) that after twisting a toric structure on
  $E$ by a character $\chi$ both sets of characters $H^1(E)$
  and $H^1(E(-1))$ are shifted simultaneously by $-\chi$ inside the lattice 
  $\L(T)$. Since $\dim H^1(E)=k-1$ and $\dim H^1(E(-1))=k$, there exists 
  a unique shift, which leads to the equality between these sums.

  So, to get the correct toric structure on a sheaf $E\bl\M^T$ we have to
  start from any toric structure on $E$, then calculate the character 
  decompositions for $H^1(E)$ and $H^1(E(-1))$, and then shift the toric 
  structure
  in order to have the equality between the sums of these characters.
  We will make all these calculations in \S4.

\SubNo{Rank 2 torsion free toric sheaves on $\PP_2$}
  The technique developed by A.K1yachko (see \cite{K1,K2}) presents a 
  description of all toric sheaves up to the isomorphisms,
  which preserve a toric structure.
  This description leads to enumeration of all connected components of $\M^T$.
  So, let us remember some combinatorial data associated with
  a toric structure on a sheaf. 

  Any sheaf $E$ on $\PP_2$ is uniquely defined by a triple of the modules 
  $E(U_i)=\Gamma(U_i;E)$. If $E$ is toric, then $T$ acts on each $E(U_i)$. 
  Under this action $t\in T$ maps a section $\MAP \s,\O,E$ to the composition
  $$\O\TO{i^{-1}_t}t_*(\O)\TO{t_*(\s)}t_*(E)\TO{\f_t}E\;,
  $$
  where $i_t$ is the tautological toric structure on $\O$ fixed above. 
  As soon as we have a torus action on
  $E(U_i)$ we can decompose $E(U_i)$ into a direct sum of eigenspaces
  $E(U_i)_\chi$ parametrized by the torus characters. Note 
  that if we multiply an eigensection of weight $\chi\in\L(T)$ by a 
  character $\xi$, which is regular on $U_i$, 
  then we get a section of weight $\chi-\xi$. So, $\fain\xi,{\s_i}$ we have
  the inclusion $E(U_i)_\chi\hookrightarrow E(U_i)_{\chi-\xi}$. 
 
\SubSubNo{\bf Triple of bifiltrations associated with a toric structure.} 
  Denote by $V$ the fiber of $E$ over the point 
  $p=(1\colon1\colon1)\in U_{012}$.
  Since any homogeneous section is uniquely determinated by
  its value at the point $p$, each $E(U_i)_\chi$ can be considered as  
  a subspace in $V$. Denote this subspace by $V^i(\chi)$. 
  As we have just seen, these subspaces form an inductive system over 
  the cone $-\s_i$,  which is opposite to $\s_i\subset\L(T)$, i.~e.
  $$V^i(\chi)\subset V^i(\chi-\xi)\qquad\fain\xi,\s_i \;.
  $$
  If we replace the character $\chi$ in the notation $V^i(\chi)$ by its 
  coordinates $m_j\bydef\sp\t_j,\chi$ and $m_k\bydef\sp\t_k,\chi$ with 
  respect to the sides of the angle $\s_i$ 
  (we always suppose that $j<k$ is the complementary to $i$ pair of indices) 
  and use the notation $V^i(m_j,m_k)$ instead of $V^i(\chi)$, 
  then we can say that for 
  $i=0,1,2$ the subspaces $V^i(m_j,m_k)$ form a triple of bifiltrations
  of $V$. These bifiltration have the properties
  $$
    V^i(m_j,m_k+1)\subset V^i(m_j,m_k)\supset V^i(m_j+1,m_k)
    \qquad\fain{m_j,m_k},{\ZZ\times\ZZ}\;,
  $$
  $$
  V^i(-\infty,-\infty)=V\qquad{\rm and}\qquad V^i(\infty,\infty)=0\;,
  $$
  and a toric structure on $E$ is uniquely defined by such 
  a triple of bifiltrations.

\SubSubNo{\bf Triple of filtrations associated with a toric structure.}
\label{filttriple} 
  Consider now the modules of sections $E(U_{ij})=\Gamma(U_i\cap U_j;E)$. 
  As above, we can associate with each of these three modules a collection of
  subspaces $V^{ij}(\chi)\subset V$ parametrized by torus characters
  $\chi\bl\L(T)$. These subspaces correspond to homogeneous sections
  of $E$ over $U_{ij}$ and satisfy the property
  $$V^{ij}(\chi-\xi)\supset V^{ij}(\chi)\qquad\fain\xi,{\s_{ij}}=
                           \{\xi\bl\L(T)\,\colon\;\sp\t_k,\xi\ge0\}\;.
  $$

  Writting $V^{ij}(m)$ instead of $V^{ij}(\chi)$, where
  $m=m_k=\sp\t_k,\chi$ (and $k$ is complementary to $(i,j)$ index),
  we get a decreasing filtration
  $$
    V=V^{ij}(-\infty)\supset\cdots\supset V^{ij}(m)\supset V^{ij}(m-1)\supset 
      \cdots\supset V^{ij}(\infty)=0
  $$
  of $V$.
  Certainly, each filtration $V^{ij}(*)$ coincides with the limits of the 
  bifiltrations  $V^i(*,*)$ and  $V^j(*,*)$:
  \begin{equation}\label{limcond}
  V^{i}(m,-\infty)=V^{j}(-\infty,m)=V^{ij}(m)\quad
                    \fain m,\ZZ\,,\;\forall\,i<j\,,\;(i,j)\subset\{0,1,2\}\;.
  \end{equation}
  Now we are in a position to formulate the following result of A.K1yachko
  extracted from \cite{K2}.

\Th\label{tsd} 
  The category of toric torsion free rank 2 sheaves with the morphisms, 
  which preserve a toric structure, is equivalent to the category of
  2-dimensional vector spaces $V$ equipped with a triple of 
  bifiltrations $V^{i}(*,*)$, $i=0,1,2$, such that any two of them have 
  the same limits (\ref{limcond}) (the morphisms in this category must preserve
  all bifiltrations).
\EF
\EP

\SubNo{Toric bundles} 
  Starting from any triple of filtrations $V^{ij}(*)$, one can construct
  a triple of bilfitrations $\overline{V^{i}}(*,*)$ defined as
  \begin{equation}\label{refcov}
  \overline{V^{i}}(m_j,m_k)=V^{ij}(m_1)\cap V^{ik}(m_2)\;.
  \end{equation}
  Evidently, this triple satisfy the conditions of \ref{tsd}. Hence, such a 
  triple always define a toric sheaf. It is easy to see that this sheaf is 
  locally free. So, we have the following result (see \cite{K1}).
 
  \Th
  The category of all rank 2 vector bundles with the
  morphisms, which preserve the toric structure, is equivalent to the 
  category of 2-dimensional vector spaces equipped with a triple of 
  filtrations $V^{ij}(*)$, $i<j$, $(i,j)\subset\{0,1,2\}$ (the morphisms
  in this category must preserve all the filtrations).
  \EF
  \EP
  
  \noindent
  Since $\dim V=2$, each of the filtration $V^{jk}(*)$ can be given by
  the following data:

  \smallskip\begingroup\leftskip=1cm\par
  \noindent\llap{--}
  1-dimensional subspace $L^i\subset V$;
  
  \smallskip\par\noindent\llap{--}
  a number $s_i\bl\ZZ$ such that 
      $V^{jk}(s_i)=V\qquad{\rm and}\qquad V^{jk}(s_i+1)=L^i\;;$
  
  \smallskip\par\noindent\llap{--}
   a number $a_i\bl\ZZ$, $a_i\ge0$ such that 
   $V^{jk}(s_i+a_i)=L^i\qquad{\rm and}\qquad V^{jk}(s_i+a_i+1)=0$
   (in other words $a_i$ is the number of 1-dimensional terms of the 
   filtration).

  \smallskip\endgroup

  \noindent
  In terms of this data the action of isomorphisms, which preserve a toric 
  structure, coincides with the action of $\PGL2,V$ on a triple of lines $L^i$.
  There is also the very nice stability criterion (see \cite{K1,K2}):

  \Th
  Locally free $\mu$-stable rank 2 toric sheaf $E$ is $\mu$-stable if and only 
  if the corresponding 1-dimensional subspaces $L^i$ are pairwise different and
  $a_i$ are the side lengths of a triangle (i.~e. they are positive integers
  satisfying the triangle inequalities). 
  \EF
  \EP

  Since $\PGL3,V$ acts transitively on the triples of pairwise 
  different 1-dimensional subspaces $L^i$, a stable toric bundle does not depend 
  (up to isomorphism preserving a toric structure) on the choose of
  the subspaces $L^i$ in $V$. So, a stable toric bundle $E$ is defined in fact
  only by six numbers $a_i$, $s_i$.  For example, the Chern classes of
  $E$ can be recovered from these numbers by the formulas (see \cite{K1})
  \begin{eqnarray*}
          c_1(E) &=& 2(s_0+s_1+s_2)\, +\, a_0+a_1+a_2\\
          c_1(E)^2-4c_2(E) &=& (-a_0+a_1+a_2)^2-4a_1a_2
  \end{eqnarray*}

  Moreover, a toric structure on $E$ depends on (and only on!) the numbers $s_i$.
  Namely, if a toric structure is twisted by a character $\chi$, then
  the numbers $a_i$ remain to be the same and the numbers $s_i$ are
  shifted by the rule
  $$s_i\longmapsto s_i-\sp\t_i,\chi
  $$
  (so, we see that the sum $s_0+s_1+s_2$ is not changed too).   
  
  We get

  \Cl\label{asdata}
  There is a 1-1 correspondence between the stable toric bundles 
  (considered up to toric isomorphisms) with 
  $(\rk,c_1,c_2)=(2,-1,k)$ and the collections
  $$(s_1,s_2,a_0,a_1,a_2)\subset\ZZ\times\ZZ\times\NN\times\NN\times\NN
  $$
  such that $a_0$, $a_1$, $a_2$ satisfy the triangle inequalities 
  and the equality
  $$ (-a_0+a_1+a_2)^2-4a_1a_2=1-4k\;.
  $$  
  If a toric structure on a bundle is changed by twisting by a character 
  $\chi\bl\L(T)$, then the triple $(a_0,a_1,a_2)$ remains to be the same
  and a pair $(s_1,s_2)$ is shifted by the vector with coordinates
  $(-\sp\t_1,\chi,-\sp\t_2,\chi)$.
  \EF
  \EP

\SubNo{Non-reflexive toric sheaves}   
  From the canonical exact triple
  $$0\TO{}E\TO{}E^{**}\TO{}C_E\TO{}0
  $$
  it follows that if $E$ is a stable toric sheaf with $(\rk,c_1,c_2)=(2,-1,k)$, 
  then $E^{**}$ is a stable toric bundle with $(\rk,c_1,c_2)=(2,-1,k-d)$,
  where $d=\dim H^0(C_E)$. The cokernel $C_E=E/E^{**}$ splits 
  into direct sum of three torsion toric sheaves $C_i$ with
  ${\rm Supp}(C_i)=p_i$, where  $p_0=(1\colon0\colon0)$, $p_1=(0\colon1\colon0)$, 
  and $p_2=(0\colon0\colon1)$ are three fixed points for the torus action
  on $\PP_2$.  

  Hence, the combinatorial description of an
  arbitrary toric sheaf consists of a non-negative integer $d$, which
  gives  ``a jump'' $c_2(E)-c_2(E^{**})$, five integers 
  $(s_1,s_2,a_0,a_1,a_2)$, which give a bundle $E^{**}$ with $c_2(E^{**})=k-d$,
  and a combinatorial description of three  $\O_{\PP_2,p_i}$-modules 
  $C_i$. 

\SubSubNo{\bf Combinatorial data and moduli of $C_i$.}
  Suppose $E$ to be given by a triple of bifiltrations $V^{i}(*,*)$. 
  Then  construct the
  filtrations $V^{ij}(*)$ defined by (\ref{limcond}). Then make from them
  the new bifiltrations $\overline{V^{i}}(*,*)$ defined by (\ref{refcov}). 
  It is not difficult to see that the vector bundle, which corresponds
  to this final bifiltration's triple, coincides with $E^{**}$. Moreover, 
  the canonical inclusion
  $E\hookrightarrow E^{**}$ corresponds to the natural embedding
  $$
  V^{i}(*,*)\subset\overline{V^{i}}(*,*)=
  V^{i}(-\infty,*)\cap V^{i}(*,-\infty)\;,
  $$
  Hence,  each $C_i$ can be decomposed via torus action by the formula
  $$ 
  C_i=\bigoplus_{m_j,m_k}
                      \overline{V^{i}}(m_j,m_k)/V^{i}(m_j,m_k)\;.
  $$
  
  This decomposition can be represented by the following picture.
  Consider an infinite table with cells indexed by the pairs $(m_j,m_k)$
  and put in each cell the corresponding space $\overline{V^i}(m_j,m_k)$.
  Such a table is divided into several natural zones, which contain the
  spaces of the same dimension. This looks like

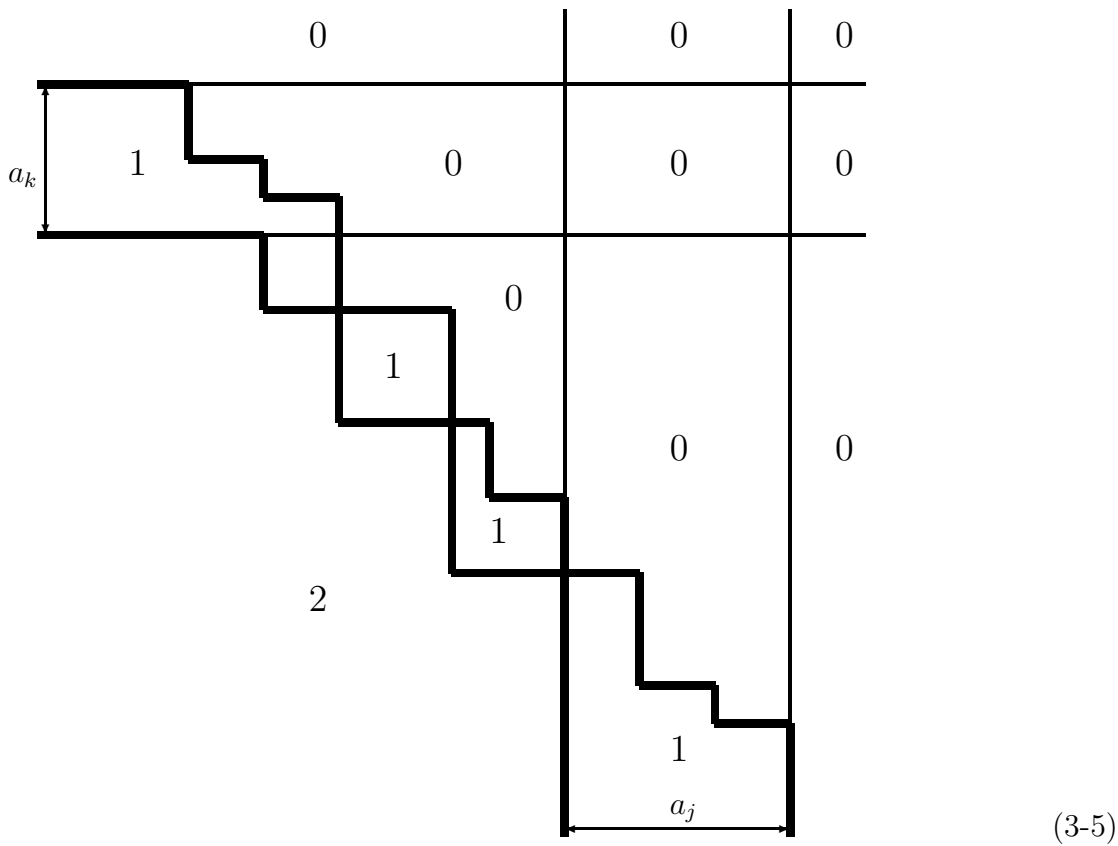
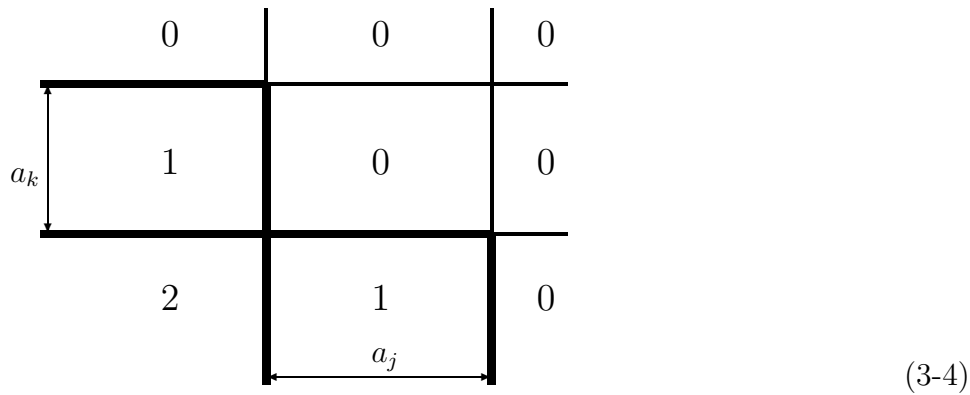
\begin{equation}\label{refbfpic}
\mbox{
\begin{picture}(70,50)
     \put(16,45){\large 0}
     \put(44,45){\large 0}
     \put(66,45){\large 0}
     \put(16,28){\large 1}
     \put(44,28){\large 0}
     \put(66,28){\large 0}
     \put(16,10){\large 2}
     \put(44,10){\large 1}
     \put(66,10){\large 0}
     \put(1,30){\vector(0,1){10}}
     \put(1,30){\vector(0,-1){10}}
     \put(45,1){\vector(1,0){15}}
     \put(45,1){\vector(-1,0){15}}
     \put(-4,27){$a_k$}
     \put(44,3){$a_j$}
     \linethickness{.4mm}
     \put(0,20){\line(1,0){70}}
     \put(0,40){\line(1,0){70}}
     \put(30,0){\line(0,1){50}}
     \put(60,0){\line(0,1){50}}
     \linethickness{1mm}

     \put(0,20){\line(1,0){60}}
     \put(0,40){\line(1,0){30}}
     \put(30,0){\line(0,1){40}}
     \put(60,0){\line(0,1){20}}
\end{picture}
}
\end{equation}

\noindent
 where the numbers indicate the dimensions of the corresponding spaces. 
 Note that the widths of 1-dimensional zones
 coincide with the numbers $a_j$, $a_k$, which define the bundle $E^{**}$
 via \ref{asdata}. The numbers $s_i$ define the placement of the
 picture: upper cells of the horizontal 1-dimensional strip have $m_k=s_k+a_k$
 and the right cells of the vertical 1-dimensional zone have $m_j=s_j+a_j$.

\begin{figure}
\begin{equation}\label{bfpic}
\mbox{
\begin{picture}(110,110)
     \put(36,105){\large 0}
     \put(84,105){\large 0}
     \put(106,105){\large 0}
     \put(62,70){\large 0}
     \put(12,88){\large 1}
     \put(54,88){\large 0}
     \put(84,88){\large 0}
     \put(106,88){\large 0}
     \put(46,61){\large 1}
     \put(60,39){\large 1}
     \put(84,50){\large 0}
     \put(106,50){\large 0}
     \put(36,30){\large 2}
     \put(84,10){\large 1}
     \put(1,90){\vector(0,1){10}}
     \put(1,90){\vector(0,-1){10}}
     \put(85,1){\vector(1,0){15}}
     \put(85,1){\vector(-1,0){15}}
     \put(-4,87){$a_k$}
     \put(84,3){$a_j$}
     \linethickness{.4mm}
     \put(0,100){\line(1,0){110}}
     \put(0,80){\line(1,0){110}}
     \put(70,0){\line(0,1){110}}
     \put(100,0){\line(0,1){110}}
     \linethickness{1mm}
     \put(0,100){\line(1,0){20}}
     \put(20,90){\line(1,0){10}}
     \put(30,85){\line(1,0){10}}
     \put(0,80){\line(1,0){30}}
     \put(30,70){\line(1,0){25}}
     \put(40,55){\line(1,0){20}}
     \put(60,45){\line(1,0){10}}
     \put(55,35){\line(1,0){25}}
     \put(80,20){\line(1,0){10}}
     \put(90,15){\line(1,0){10}}
     \put(20,90){\line(0,1){10}}
     \put(30,70){\line(0,1){10}}
     \put(30,85){\line(0,1){5}}
     \put(40,55){\line(0,1){30}}
     \put(55,35){\line(0,1){35}}
     \put(70,0){\line(0,1){45}}
     \put(80,20){\line(0,1){15}}
     \put(90,15){\line(0,1){5}}
     \put(100,0){\line(0,1){15}}
     \put(60,45){\line(0,1){10}}
\end{picture}
}
\end{equation}\end{figure}

  Similar table for the bifiltration $V^i(*,*)$ is obtained from the 
  previous one by changing some dimensions by smaler like in
  (\ref{bfpic}). An important property of the 1-dimensional zone of this 
  table is that any two cells, which have a common side, have to contain 
  the same 1-dimensional spaces. So, 1-dimensional zone splits into 
  {\it maximal connected components\/}: two cells 
  belong to the same component iff they can be connected inside this 
  component by a sequence of cells such that any two consequent elements 
  of this sequence have a common side. All cells of each connected 
  component contain the same 1-dimensional subspace.
  
  Note, that in any case there are exactly two not bounded components. 
  They contain 1-dimensional  spaces, which come from the previous picture 
  (i.e. from $E^{**}$ -- corresponding cells form the 
  bottom end of the vertical strip and the left end of the horizontal strip). 
  These two 1-dimensional spaces must be different, because of stability, 
  and we can consider them to be fixed by some toric isomorphism. 

  Other connected components are bounded (there are 2 such components 
  in (\ref{bfpic})). There are no
  restrictions on the 1-dimensional spaces placed in each of these 
  components. Hence, each torsion sheaf $C_i$ defines and is uniquely 
  defined by a table like (\ref{bfpic}) and a point from 
  $$\underbrace{\PP_1\times\PP_1\times\cdots\times\PP_1}_{N_i}\;, 
  $$ 
  where $N_i$ is the number of bounded connected components in 
  1-dimensional zone and each $\PP_1=\PP(V)$ parametrize the choice of 
  1-dimensional subspace in $V$ placed in the corresponding connected 
  component.

  So, each connected component $Y\subset\M^T$ is uniquely defined by
  a triple of numbers $a_i$, which give $E^{**}$ (the same for all $E\bl Y$),
  a triple of numbers $d_i=\dim H^0(C_i)$, and a triple of pictures like
  (\ref{bfpic}). Such a component is isomorphic to the direct product of
  $N=N_0+N_1+N_2$ projective lines, which parametrize the choice of subspaces
  placed in bounded 1-dimensional zones of the tables (\ref{bfpic}).

\SubSubNo{\bf Character decomposition and eigensubbundles in $H^0(C_i)$.}
\label{hdecomp}
  The picture, which gives the character decomposition for $H^0(C_i)$ 
  follows immediately from (\ref{bfpic}), (\ref{refbfpic}): the character
  $\chi$ with the coordinates $m_j=\sp\t_j,\chi$, $m_k=\sp\t_k,\chi$ is present
  in $H^0(C_i)$ if and only if a subspace placed in the $(m_j,m_k)$-cell of 
  (\ref{bfpic}) is smaler then the one placed in the same cell of 
  (\ref{refbfpic}). The difference between
  the dimensions of these subspaces equals the multiplicity of $\chi$.

  When $E$ is runnig through a connected component $Y\subset\M^T$ the
  eigenspace corresponding to $\chi$ form a vector bundle over 
  $Y=\PP_1\times\PP_1\times\cdots\times\PP_1$. 
  This bundle is nontrivial iff $\chi$ comes from a bounded 1-dimensional 
  zone of (\ref{bfpic}). In this case it equals the pull-back of the
  universal factor line bundle $\O_{\PP_1}(1)$ over the $\PP_1$-multiplier,
  which corresponds to the bounded 1-dimensional zone of (\ref{bfpic}) 
  what $\chi$ comes from.

\SubSubNo{\bf Description of $H^0(C_i)$ by a Young diagram's pair.}
\label{ydpairs}
  Combinatorially, it is convenient to represent the
  character table for $H^0(C_i)$ as a sum of two Young diagrams 
  $\l_i$, $\mu_i$ filled by 1-dimensional spaces like in (\ref{picsum})
  ($\l_i$, $\mu_i$ are any satisfying the condition 
  $|\l_i|+|\mu_i|=d_i=\dim H^0(C_i)$). 
\begin{figure}
\begin{equation}\label{picsum}
\setlength{\unitlength}{0.075in}
\begin{picture}(68.46,41.70)
\special{em:linewidth 0.014in}
\put(64,35){{\setbox0=\hbox{$\mu_i$}\kern-\wd0\raise\dp0\box0}}
\put(10,35){{\setbox0=\hbox{$\lambda_i$}\kern-\wd0\raise\dp0\box0}}
\put(30,7){{\setbox0=\hbox{$H^0(C_i)$}\kern-\wd0\raise\dp0\box0}}
\put(28.26,26.70){\special{em:moveto}}
\put(48.26,26.70){\special{em:lineto}}
\put(48.26,5.03){\special{em:lineto}}
\put(46.66,5.03){\special{em:lineto}}
\put(46.66,10.03){\special{em:lineto}}
\put(45.00,10.03){\special{em:lineto}}
\put(44.96,13.36){\special{em:lineto}}
\put(40.00,13.36){\special{em:lineto}}
\put(40.00,16.70){\special{em:lineto}}
\put(33.33,16.70){\special{em:lineto}}
\put(33.33,19.96){\special{em:lineto}}
\put(31.66,19.96){\special{em:lineto}}
\put(31.63,25.00){\special{em:lineto}}
\put(28.26,25.00){\special{em:lineto}}
\put(28.26,26.70){\special{em:lineto}}
\put(28.26,26.70){\special{em:moveto}}
\put(48.26,26.70){\special{em:lineto}}
\put(48.26,5.03){\special{em:lineto}}
\put(46.66,5.03){\special{em:lineto}}
\put(46.66,10.03){\special{em:lineto}}
\put(45.00,10.03){\special{em:lineto}}
\put(44.96,13.36){\special{em:lineto}}
\put(40.00,13.36){\special{em:lineto}}
\put(40.00,16.70){\special{em:lineto}}
\put(33.33,16.70){\special{em:lineto}}
\put(33.33,19.96){\special{em:lineto}}
\put(31.66,19.96){\special{em:lineto}}
\put(31.63,25.00){\special{em:lineto}}
\put(28.26,25.00){\special{em:lineto}}
\put(28.26,26.70){\special{em:lineto}}
\put(33.40,21.66){\special{em:moveto}}
\put(53.40,21.66){\special{em:lineto}}
\put(53.36,3.40){\special{em:lineto}}
\put(50.06,3.40){\special{em:lineto}}
\put(50.06,5.00){\special{em:lineto}}
\put(46.70,5.00){\special{em:lineto}}
\put(46.76,6.66){\special{em:lineto}}
\put(45.06,6.70){\special{em:lineto}}
\put(45.10,10.00){\special{em:lineto}}
\put(40.03,10.00){\special{em:lineto}}
\put(40.03,13.36){\special{em:lineto}}
\put(38.40,13.36){\special{em:lineto}}
\put(38.40,16.66){\special{em:lineto}}
\put(36.73,16.66){\special{em:lineto}}
\put(36.73,20.00){\special{em:lineto}}
\put(33.40,20.00){\special{em:lineto}}
\put(33.40,21.66){\special{em:lineto}}
\put(33.40,21.66){\special{em:moveto}}
\put(53.40,21.66){\special{em:lineto}}
\put(53.36,3.40){\special{em:lineto}}
\put(50.06,3.40){\special{em:lineto}}
\put(50.06,5.00){\special{em:lineto}}
\put(46.70,5.00){\special{em:lineto}}
\put(46.76,6.66){\special{em:lineto}}
\put(45.06,6.70){\special{em:lineto}}
\put(45.10,10.00){\special{em:lineto}}
\put(40.03,10.00){\special{em:lineto}}
\put(40.03,13.36){\special{em:lineto}}
\put(38.40,13.36){\special{em:lineto}}
\put(38.40,16.66){\special{em:lineto}}
\put(36.73,16.66){\special{em:lineto}}
\put(36.73,20.00){\special{em:lineto}}
\put(33.40,20.00){\special{em:lineto}}
\put(33.40,21.66){\special{em:lineto}}
\put(0.00,41.70){\special{em:moveto}}
\put(20.00,41.70){\special{em:lineto}}
\put(20.00,20.03){\special{em:lineto}}
\put(18.40,20.03){\special{em:lineto}}
\put(18.40,25.03){\special{em:lineto}}
\put(16.73,25.03){\special{em:lineto}}
\put(16.70,28.36){\special{em:lineto}}
\put(11.73,28.36){\special{em:lineto}}
\put(11.73,31.70){\special{em:lineto}}
\put(5.06,31.70){\special{em:lineto}}
\put(5.06,34.96){\special{em:lineto}}
\put(3.40,34.96){\special{em:lineto}}
\put(3.36,40.00){\special{em:lineto}}
\put(0.00,40.00){\special{em:lineto}}
\put(0.00,41.70){\special{em:lineto}}
\put(0.00,41.70){\special{em:moveto}}
\put(20.00,41.70){\special{em:lineto}}
\put(20.00,20.03){\special{em:lineto}}
\put(18.40,20.03){\special{em:lineto}}
\put(18.40,25.03){\special{em:lineto}}
\put(16.73,25.03){\special{em:lineto}}
\put(16.70,28.36){\special{em:lineto}}
\put(11.73,28.36){\special{em:lineto}}
\put(11.73,31.70){\special{em:lineto}}
\put(5.06,31.70){\special{em:lineto}}
\put(5.06,34.96){\special{em:lineto}}
\put(3.40,34.96){\special{em:lineto}}
\put(3.36,40.00){\special{em:lineto}}
\put(0.00,40.00){\special{em:lineto}}
\put(0.00,41.70){\special{em:lineto}}
\put(48.46,41.50){\special{em:moveto}}
\put(68.46,41.50){\special{em:lineto}}
\put(68.43,23.23){\special{em:lineto}}
\put(65.13,23.23){\special{em:lineto}}
\put(65.13,24.83){\special{em:lineto}}
\put(61.76,24.83){\special{em:lineto}}
\put(61.83,26.50){\special{em:lineto}}
\put(60.13,26.53){\special{em:lineto}}
\put(60.16,29.83){\special{em:lineto}}
\put(55.10,29.83){\special{em:lineto}}
\put(55.10,33.20){\special{em:lineto}}
\put(53.46,33.20){\special{em:lineto}}
\put(53.46,36.50){\special{em:lineto}}
\put(51.80,36.50){\special{em:lineto}}
\put(51.80,39.83){\special{em:lineto}}
\put(48.46,39.83){\special{em:lineto}}
\put(48.46,41.50){\special{em:lineto}}
\put(48.46,41.50){\special{em:moveto}}
\put(68.46,41.50){\special{em:lineto}}
\put(68.43,23.23){\special{em:lineto}}
\put(65.13,23.23){\special{em:lineto}}
\put(65.13,24.83){\special{em:lineto}}
\put(61.76,24.83){\special{em:lineto}}
\put(61.83,26.50){\special{em:lineto}}
\put(60.13,26.53){\special{em:lineto}}
\put(60.16,29.83){\special{em:lineto}}
\put(55.10,29.83){\special{em:lineto}}
\put(55.10,33.20){\special{em:lineto}}
\put(53.46,33.20){\special{em:lineto}}
\put(53.46,36.50){\special{em:lineto}}
\put(51.80,36.50){\special{em:lineto}}
\put(51.80,39.83){\special{em:lineto}}
\put(48.46,39.83){\special{em:lineto}}
\put(48.46,41.50){\special{em:lineto}}
\put(33.46,21.63){\special{em:moveto}}
\put(33.36,19.96){\special{em:lineto}}
\put(36.73,19.96){\special{em:lineto}}
\put(36.70,16.66){\special{em:lineto}}
\put(39.93,16.73){\special{em:lineto}}
\put(40.03,13.33){\special{em:lineto}}
\put(44.93,13.33){\special{em:lineto}}
\put(44.93,10.00){\special{em:lineto}}
\put(46.56,10.00){\special{em:lineto}}
\put(46.66,5.03){\special{em:lineto}}
\put(48.40,5.03){\special{em:lineto}}
\put(48.26,21.63){\special{em:lineto}}
\put(33.46,21.63){\special{em:lineto}}
\put(33.46,21.63){\special{em:moveto}}
\put(33.36,19.96){\special{em:lineto}}
\put(36.73,19.96){\special{em:lineto}}
\put(36.70,16.66){\special{em:lineto}}
\put(39.93,16.73){\special{em:lineto}}
\put(40.03,13.33){\special{em:lineto}}
\put(44.93,13.33){\special{em:lineto}}
\put(44.93,10.00){\special{em:lineto}}
\put(46.56,10.00){\special{em:lineto}}
\put(46.66,5.03){\special{em:lineto}}
\put(48.40,5.03){\special{em:lineto}}
\put(48.26,21.63){\special{em:lineto}}
\put(33.46,21.63){\special{em:lineto}}
\put(53.43,26.73){\special{em:moveto}}
\put(53.43,26.65){\special{em:lineto}}
\put(53.43,26.40){\special{em:moveto}}
\put(53.43,26.31){\special{em:lineto}}
\put(53.43,26.06){\special{em:moveto}}
\put(53.43,25.98){\special{em:lineto}}
\put(53.43,25.73){\special{em:moveto}}
\put(53.43,25.65){\special{em:lineto}}
\put(53.43,25.40){\special{em:moveto}}
\put(53.43,25.31){\special{em:lineto}}
\put(53.43,25.06){\special{em:moveto}}
\put(53.43,24.98){\special{em:lineto}}
\put(53.43,24.73){\special{em:moveto}}
\put(53.43,24.65){\special{em:lineto}}
\put(53.43,24.40){\special{em:moveto}}
\put(53.43,24.31){\special{em:lineto}}
\put(53.43,24.06){\special{em:moveto}}
\put(53.43,23.98){\special{em:lineto}}
\put(53.43,23.73){\special{em:moveto}}
\put(53.43,23.65){\special{em:lineto}}
\put(53.43,23.40){\special{em:moveto}}
\put(53.43,23.31){\special{em:lineto}}
\put(53.43,23.06){\special{em:moveto}}
\put(53.43,22.98){\special{em:lineto}}
\put(53.43,22.73){\special{em:moveto}}
\put(53.43,22.65){\special{em:lineto}}
\put(53.43,22.40){\special{em:moveto}}
\put(53.43,22.31){\special{em:lineto}}
\put(53.43,22.06){\special{em:moveto}}
\put(53.43,21.98){\special{em:lineto}}
\put(53.43,21.73){\special{em:moveto}}
\put(53.43,21.65){\special{em:lineto}}
\put(53.43,21.40){\special{em:moveto}}
\put(53.43,21.31){\special{em:lineto}}
\put(53.43,21.06){\special{em:moveto}}
\put(53.43,20.98){\special{em:lineto}}
\put(53.43,20.73){\special{em:moveto}}
\put(53.43,20.65){\special{em:lineto}}
\put(53.43,20.40){\special{em:moveto}}
\put(53.43,20.31){\special{em:lineto}}
\put(53.43,20.06){\special{em:moveto}}
\put(53.43,19.98){\special{em:lineto}}
\put(53.43,19.73){\special{em:moveto}}
\put(53.43,19.65){\special{em:lineto}}
\put(53.43,19.40){\special{em:moveto}}
\put(53.43,19.31){\special{em:lineto}}
\put(53.43,19.06){\special{em:moveto}}
\put(53.43,18.98){\special{em:lineto}}
\put(53.43,18.73){\special{em:moveto}}
\put(53.43,18.65){\special{em:lineto}}
\put(53.43,18.40){\special{em:moveto}}
\put(53.43,18.31){\special{em:lineto}}
\put(53.43,18.06){\special{em:moveto}}
\put(53.43,17.98){\special{em:lineto}}
\put(53.43,17.73){\special{em:moveto}}
\put(53.43,17.65){\special{em:lineto}}
\put(53.43,17.40){\special{em:moveto}}
\put(53.43,17.31){\special{em:lineto}}
\put(53.43,17.06){\special{em:moveto}}
\put(53.43,16.98){\special{em:lineto}}
\put(53.43,16.73){\special{em:moveto}}
\put(53.43,16.65){\special{em:lineto}}
\put(53.43,16.40){\special{em:moveto}}
\put(53.43,16.31){\special{em:lineto}}
\put(53.43,16.06){\special{em:moveto}}
\put(53.43,15.98){\special{em:lineto}}
\put(53.43,15.73){\special{em:moveto}}
\put(53.43,15.65){\special{em:lineto}}
\put(53.43,15.40){\special{em:moveto}}
\put(53.43,15.31){\special{em:lineto}}
\put(53.43,15.06){\special{em:moveto}}
\put(53.43,14.98){\special{em:lineto}}
\put(53.43,14.73){\special{em:moveto}}
\put(53.43,14.65){\special{em:lineto}}
\put(53.43,14.40){\special{em:moveto}}
\put(53.43,14.31){\special{em:lineto}}
\put(53.43,14.06){\special{em:moveto}}
\put(53.43,13.98){\special{em:lineto}}
\put(53.43,13.73){\special{em:moveto}}
\put(53.43,13.65){\special{em:lineto}}
\put(53.43,13.40){\special{em:moveto}}
\put(53.43,13.31){\special{em:lineto}}
\put(53.43,13.06){\special{em:moveto}}
\put(53.43,12.98){\special{em:lineto}}
\put(53.43,12.73){\special{em:moveto}}
\put(53.43,12.65){\special{em:lineto}}
\put(53.43,12.40){\special{em:moveto}}
\put(53.43,12.31){\special{em:lineto}}
\put(53.43,12.06){\special{em:moveto}}
\put(53.43,11.98){\special{em:lineto}}
\put(53.43,11.73){\special{em:moveto}}
\put(53.43,11.65){\special{em:lineto}}
\put(53.43,11.40){\special{em:moveto}}
\put(53.43,11.31){\special{em:lineto}}
\put(53.43,11.06){\special{em:moveto}}
\put(53.43,10.98){\special{em:lineto}}
\put(53.43,10.73){\special{em:moveto}}
\put(53.43,10.65){\special{em:lineto}}
\put(53.43,10.40){\special{em:moveto}}
\put(53.43,10.31){\special{em:lineto}}
\put(53.43,10.06){\special{em:moveto}}
\put(53.43,9.98){\special{em:lineto}}
\put(53.43,9.73){\special{em:moveto}}
\put(53.43,9.65){\special{em:lineto}}
\put(53.43,9.40){\special{em:moveto}}
\put(53.43,9.31){\special{em:lineto}}
\put(53.43,9.06){\special{em:moveto}}
\put(53.43,8.98){\special{em:lineto}}
\put(53.43,8.73){\special{em:moveto}}
\put(53.43,8.65){\special{em:lineto}}
\put(53.43,8.40){\special{em:moveto}}
\put(53.43,8.31){\special{em:lineto}}
\put(53.43,8.06){\special{em:moveto}}
\put(53.43,7.98){\special{em:lineto}}
\put(53.43,7.73){\special{em:moveto}}
\put(53.43,7.65){\special{em:lineto}}
\put(53.43,7.40){\special{em:moveto}}
\put(53.43,7.31){\special{em:lineto}}
\put(53.43,7.06){\special{em:moveto}}
\put(53.43,6.98){\special{em:lineto}}
\put(53.43,6.73){\special{em:moveto}}
\put(53.43,6.65){\special{em:lineto}}
\put(53.43,6.40){\special{em:moveto}}
\put(53.43,6.31){\special{em:lineto}}
\put(53.43,6.06){\special{em:moveto}}
\put(53.43,5.98){\special{em:lineto}}
\put(53.43,5.73){\special{em:moveto}}
\put(53.43,5.65){\special{em:lineto}}
\put(53.43,5.40){\special{em:moveto}}
\put(53.43,5.31){\special{em:lineto}}
\put(53.43,5.06){\special{em:moveto}}
\put(53.43,4.98){\special{em:lineto}}
\put(53.43,4.73){\special{em:moveto}}
\put(53.43,4.65){\special{em:lineto}}
\put(53.43,4.40){\special{em:moveto}}
\put(53.43,4.31){\special{em:lineto}}
\put(53.43,4.06){\special{em:moveto}}
\put(53.43,3.98){\special{em:lineto}}
\put(53.43,3.73){\special{em:moveto}}
\put(53.43,3.65){\special{em:lineto}}
\put(53.43,3.40){\special{em:moveto}}
\put(53.43,3.31){\special{em:lineto}}
\put(53.43,3.06){\special{em:moveto}}
\put(53.43,2.98){\special{em:lineto}}
\put(53.43,2.73){\special{em:moveto}}
\put(53.43,2.65){\special{em:lineto}}
\put(53.43,2.40){\special{em:moveto}}
\put(53.43,2.31){\special{em:lineto}}
\put(53.43,2.06){\special{em:moveto}}
\put(53.43,1.98){\special{em:lineto}}
\put(53.43,1.73){\special{em:moveto}}
\put(53.43,1.65){\special{em:lineto}}
\put(53.43,1.40){\special{em:moveto}}
\put(53.43,1.31){\special{em:lineto}}
\put(53.43,1.06){\special{em:moveto}}
\put(53.43,0.98){\special{em:lineto}}
\put(53.41,0.73){\special{em:moveto}}
\put(53.41,0.65){\special{em:lineto}}
\put(53.40,0.40){\special{em:moveto}}
\put(53.40,0.31){\special{em:lineto}}
\put(53.38,0.06){\special{em:moveto}}
\put(53.36,0.03){\special{em:lineto}}
\put(48.40,26.66){\special{em:moveto}}
\put(48.40,26.58){\special{em:lineto}}
\put(48.40,26.33){\special{em:moveto}}
\put(48.40,26.25){\special{em:lineto}}
\put(48.40,26.00){\special{em:moveto}}
\put(48.40,25.91){\special{em:lineto}}
\put(48.40,25.66){\special{em:moveto}}
\put(48.40,25.58){\special{em:lineto}}
\put(48.40,25.33){\special{em:moveto}}
\put(48.40,25.25){\special{em:lineto}}
\put(48.40,25.00){\special{em:moveto}}
\put(48.40,24.91){\special{em:lineto}}
\put(48.40,24.66){\special{em:moveto}}
\put(48.40,24.58){\special{em:lineto}}
\put(48.40,24.33){\special{em:moveto}}
\put(48.40,24.25){\special{em:lineto}}
\put(48.40,24.00){\special{em:moveto}}
\put(48.40,23.91){\special{em:lineto}}
\put(48.40,23.66){\special{em:moveto}}
\put(48.40,23.58){\special{em:lineto}}
\put(48.40,23.33){\special{em:moveto}}
\put(48.40,23.25){\special{em:lineto}}
\put(48.40,23.00){\special{em:moveto}}
\put(48.40,22.91){\special{em:lineto}}
\put(48.40,22.66){\special{em:moveto}}
\put(48.40,22.58){\special{em:lineto}}
\put(48.40,22.33){\special{em:moveto}}
\put(48.40,22.25){\special{em:lineto}}
\put(48.40,22.00){\special{em:moveto}}
\put(48.40,21.91){\special{em:lineto}}
\put(48.40,21.66){\special{em:moveto}}
\put(48.40,21.58){\special{em:lineto}}
\put(48.40,21.33){\special{em:moveto}}
\put(48.40,21.25){\special{em:lineto}}
\put(48.40,21.00){\special{em:moveto}}
\put(48.40,20.91){\special{em:lineto}}
\put(48.40,20.66){\special{em:moveto}}
\put(48.40,20.58){\special{em:lineto}}
\put(48.40,20.33){\special{em:moveto}}
\put(48.40,20.25){\special{em:lineto}}
\put(48.40,20.00){\special{em:moveto}}
\put(48.40,19.91){\special{em:lineto}}
\put(48.40,19.66){\special{em:moveto}}
\put(48.40,19.58){\special{em:lineto}}
\put(48.40,19.33){\special{em:moveto}}
\put(48.40,19.25){\special{em:lineto}}
\put(48.40,19.00){\special{em:moveto}}
\put(48.40,18.91){\special{em:lineto}}
\put(48.40,18.66){\special{em:moveto}}
\put(48.40,18.58){\special{em:lineto}}
\put(48.40,18.33){\special{em:moveto}}
\put(48.40,18.25){\special{em:lineto}}
\put(48.40,18.00){\special{em:moveto}}
\put(48.40,17.91){\special{em:lineto}}
\put(48.40,17.66){\special{em:moveto}}
\put(48.40,17.58){\special{em:lineto}}
\put(48.40,17.33){\special{em:moveto}}
\put(48.40,17.25){\special{em:lineto}}
\put(48.40,17.00){\special{em:moveto}}
\put(48.40,16.91){\special{em:lineto}}
\put(48.40,16.66){\special{em:moveto}}
\put(48.40,16.58){\special{em:lineto}}
\put(48.40,16.33){\special{em:moveto}}
\put(48.40,16.25){\special{em:lineto}}
\put(48.40,16.00){\special{em:moveto}}
\put(48.40,15.91){\special{em:lineto}}
\put(48.40,15.66){\special{em:moveto}}
\put(48.40,15.58){\special{em:lineto}}
\put(48.40,15.33){\special{em:moveto}}
\put(48.40,15.25){\special{em:lineto}}
\put(48.40,15.00){\special{em:moveto}}
\put(48.40,14.91){\special{em:lineto}}
\put(48.40,14.66){\special{em:moveto}}
\put(48.40,14.58){\special{em:lineto}}
\put(48.40,14.33){\special{em:moveto}}
\put(48.40,14.25){\special{em:lineto}}
\put(48.40,14.00){\special{em:moveto}}
\put(48.40,13.91){\special{em:lineto}}
\put(48.40,13.66){\special{em:moveto}}
\put(48.40,13.58){\special{em:lineto}}
\put(48.40,13.33){\special{em:moveto}}
\put(48.40,13.25){\special{em:lineto}}
\put(48.40,13.00){\special{em:moveto}}
\put(48.40,12.91){\special{em:lineto}}
\put(48.40,12.66){\special{em:moveto}}
\put(48.40,12.58){\special{em:lineto}}
\put(48.40,12.33){\special{em:moveto}}
\put(48.40,12.25){\special{em:lineto}}
\put(48.40,12.00){\special{em:moveto}}
\put(48.40,11.91){\special{em:lineto}}
\put(48.40,11.66){\special{em:moveto}}
\put(48.40,11.58){\special{em:lineto}}
\put(48.40,11.33){\special{em:moveto}}
\put(48.40,11.25){\special{em:lineto}}
\put(48.40,11.00){\special{em:moveto}}
\put(48.40,10.91){\special{em:lineto}}
\put(48.40,10.66){\special{em:moveto}}
\put(48.40,10.58){\special{em:lineto}}
\put(48.40,10.33){\special{em:moveto}}
\put(48.40,10.25){\special{em:lineto}}
\put(48.40,10.00){\special{em:moveto}}
\put(48.40,9.91){\special{em:lineto}}
\put(48.40,9.66){\special{em:moveto}}
\put(48.40,9.58){\special{em:lineto}}
\put(48.40,9.33){\special{em:moveto}}
\put(48.40,9.25){\special{em:lineto}}
\put(48.40,9.00){\special{em:moveto}}
\put(48.40,8.91){\special{em:lineto}}
\put(48.40,8.66){\special{em:moveto}}
\put(48.40,8.58){\special{em:lineto}}
\put(48.40,8.33){\special{em:moveto}}
\put(48.40,8.25){\special{em:lineto}}
\put(48.40,8.00){\special{em:moveto}}
\put(48.40,7.91){\special{em:lineto}}
\put(48.40,7.66){\special{em:moveto}}
\put(48.40,7.58){\special{em:lineto}}
\put(48.40,7.33){\special{em:moveto}}
\put(48.40,7.25){\special{em:lineto}}
\put(48.40,7.00){\special{em:moveto}}
\put(48.40,6.91){\special{em:lineto}}
\put(48.40,6.66){\special{em:moveto}}
\put(48.40,6.58){\special{em:lineto}}
\put(48.40,6.33){\special{em:moveto}}
\put(48.40,6.25){\special{em:lineto}}
\put(48.40,6.00){\special{em:moveto}}
\put(48.40,5.91){\special{em:lineto}}
\put(48.40,5.66){\special{em:moveto}}
\put(48.40,5.58){\special{em:lineto}}
\put(48.40,5.33){\special{em:moveto}}
\put(48.40,5.25){\special{em:lineto}}
\put(48.40,5.00){\special{em:moveto}}
\put(48.40,4.91){\special{em:lineto}}
\put(48.40,4.66){\special{em:moveto}}
\put(48.40,4.58){\special{em:lineto}}
\put(48.40,4.33){\special{em:moveto}}
\put(48.40,4.25){\special{em:lineto}}
\put(48.40,4.00){\special{em:moveto}}
\put(48.40,3.91){\special{em:lineto}}
\put(48.40,3.66){\special{em:moveto}}
\put(48.40,3.58){\special{em:lineto}}
\put(48.40,3.33){\special{em:moveto}}
\put(48.40,3.25){\special{em:lineto}}
\put(48.40,3.00){\special{em:moveto}}
\put(48.40,2.91){\special{em:lineto}}
\put(48.40,2.66){\special{em:moveto}}
\put(48.40,2.58){\special{em:lineto}}
\put(48.40,2.33){\special{em:moveto}}
\put(48.40,2.25){\special{em:lineto}}
\put(48.40,2.00){\special{em:moveto}}
\put(48.40,1.91){\special{em:lineto}}
\put(48.40,1.66){\special{em:moveto}}
\put(48.40,1.58){\special{em:lineto}}
\put(48.40,1.33){\special{em:moveto}}
\put(48.40,1.25){\special{em:lineto}}
\put(48.40,1.00){\special{em:moveto}}
\put(48.40,0.91){\special{em:lineto}}
\put(48.40,0.66){\special{em:moveto}}
\put(48.40,0.58){\special{em:lineto}}
\put(48.40,0.33){\special{em:moveto}}
\put(48.40,0.25){\special{em:lineto}}
\put(23.40,21.70){\special{em:moveto}}
\put(23.48,21.70){\special{em:lineto}}
\put(23.73,21.70){\special{em:moveto}}
\put(23.81,21.70){\special{em:lineto}}
\put(24.06,21.70){\special{em:moveto}}
\put(24.15,21.70){\special{em:lineto}}
\put(24.40,21.70){\special{em:moveto}}
\put(24.48,21.70){\special{em:lineto}}
\put(24.73,21.70){\special{em:moveto}}
\put(24.81,21.70){\special{em:lineto}}
\put(25.06,21.70){\special{em:moveto}}
\put(25.15,21.70){\special{em:lineto}}
\put(25.40,21.70){\special{em:moveto}}
\put(25.48,21.70){\special{em:lineto}}
\put(25.73,21.70){\special{em:moveto}}
\put(25.81,21.70){\special{em:lineto}}
\put(26.06,21.70){\special{em:moveto}}
\put(26.15,21.70){\special{em:lineto}}
\put(26.40,21.70){\special{em:moveto}}
\put(26.48,21.70){\special{em:lineto}}
\put(26.73,21.70){\special{em:moveto}}
\put(26.81,21.70){\special{em:lineto}}
\put(27.06,21.70){\special{em:moveto}}
\put(27.15,21.70){\special{em:lineto}}
\put(27.40,21.70){\special{em:moveto}}
\put(27.48,21.70){\special{em:lineto}}
\put(27.73,21.70){\special{em:moveto}}
\put(27.81,21.70){\special{em:lineto}}
\put(28.06,21.70){\special{em:moveto}}
\put(28.15,21.70){\special{em:lineto}}
\put(28.40,21.70){\special{em:moveto}}
\put(28.48,21.70){\special{em:lineto}}
\put(28.73,21.70){\special{em:moveto}}
\put(28.81,21.70){\special{em:lineto}}
\put(29.06,21.70){\special{em:moveto}}
\put(29.15,21.70){\special{em:lineto}}
\put(29.40,21.70){\special{em:moveto}}
\put(29.48,21.70){\special{em:lineto}}
\put(29.73,21.70){\special{em:moveto}}
\put(29.81,21.70){\special{em:lineto}}
\put(30.06,21.70){\special{em:moveto}}
\put(30.15,21.70){\special{em:lineto}}
\put(30.40,21.70){\special{em:moveto}}
\put(30.48,21.70){\special{em:lineto}}
\put(30.73,21.70){\special{em:moveto}}
\put(30.81,21.70){\special{em:lineto}}
\put(31.06,21.70){\special{em:moveto}}
\put(31.15,21.70){\special{em:lineto}}
\put(31.40,21.70){\special{em:moveto}}
\put(31.48,21.70){\special{em:lineto}}
\put(31.73,21.70){\special{em:moveto}}
\put(31.81,21.70){\special{em:lineto}}
\put(32.06,21.70){\special{em:moveto}}
\put(32.15,21.70){\special{em:lineto}}
\put(32.40,21.70){\special{em:moveto}}
\put(32.48,21.70){\special{em:lineto}}
\put(32.73,21.70){\special{em:moveto}}
\put(32.81,21.70){\special{em:lineto}}
\put(33.06,21.70){\special{em:moveto}}
\put(33.15,21.70){\special{em:lineto}}
\put(33.40,21.70){\special{em:moveto}}
\put(33.48,21.70){\special{em:lineto}}
\put(33.73,21.70){\special{em:moveto}}
\put(33.81,21.70){\special{em:lineto}}
\put(34.06,21.70){\special{em:moveto}}
\put(34.15,21.70){\special{em:lineto}}
\put(34.40,21.70){\special{em:moveto}}
\put(34.48,21.70){\special{em:lineto}}
\put(34.73,21.70){\special{em:moveto}}
\put(34.81,21.70){\special{em:lineto}}
\put(35.06,21.70){\special{em:moveto}}
\put(35.15,21.70){\special{em:lineto}}
\put(35.40,21.70){\special{em:moveto}}
\put(35.48,21.70){\special{em:lineto}}
\put(35.73,21.70){\special{em:moveto}}
\put(35.81,21.70){\special{em:lineto}}
\put(36.06,21.70){\special{em:moveto}}
\put(36.15,21.70){\special{em:lineto}}
\put(36.40,21.70){\special{em:moveto}}
\put(36.48,21.70){\special{em:lineto}}
\put(36.73,21.70){\special{em:moveto}}
\put(36.81,21.70){\special{em:lineto}}
\put(37.06,21.70){\special{em:moveto}}
\put(37.15,21.70){\special{em:lineto}}
\put(37.40,21.70){\special{em:moveto}}
\put(37.48,21.70){\special{em:lineto}}
\put(37.73,21.70){\special{em:moveto}}
\put(37.81,21.70){\special{em:lineto}}
\put(38.06,21.70){\special{em:moveto}}
\put(38.15,21.70){\special{em:lineto}}
\put(38.40,21.70){\special{em:moveto}}
\put(38.48,21.70){\special{em:lineto}}
\put(38.73,21.70){\special{em:moveto}}
\put(38.81,21.70){\special{em:lineto}}
\put(39.06,21.70){\special{em:moveto}}
\put(39.15,21.70){\special{em:lineto}}
\put(39.40,21.70){\special{em:moveto}}
\put(39.48,21.70){\special{em:lineto}}
\put(39.73,21.70){\special{em:moveto}}
\put(39.81,21.70){\special{em:lineto}}
\put(40.06,21.70){\special{em:moveto}}
\put(40.15,21.70){\special{em:lineto}}
\put(40.40,21.70){\special{em:moveto}}
\put(40.48,21.70){\special{em:lineto}}
\put(40.73,21.70){\special{em:moveto}}
\put(40.81,21.70){\special{em:lineto}}
\put(41.06,21.70){\special{em:moveto}}
\put(41.15,21.70){\special{em:lineto}}
\put(41.40,21.70){\special{em:moveto}}
\put(41.48,21.70){\special{em:lineto}}
\put(41.73,21.70){\special{em:moveto}}
\put(41.81,21.70){\special{em:lineto}}
\put(42.06,21.70){\special{em:moveto}}
\put(42.15,21.70){\special{em:lineto}}
\put(42.40,21.70){\special{em:moveto}}
\put(42.48,21.70){\special{em:lineto}}
\put(42.73,21.70){\special{em:moveto}}
\put(42.81,21.70){\special{em:lineto}}
\put(43.06,21.70){\special{em:moveto}}
\put(43.15,21.70){\special{em:lineto}}
\put(43.40,21.70){\special{em:moveto}}
\put(43.48,21.70){\special{em:lineto}}
\put(43.73,21.70){\special{em:moveto}}
\put(43.81,21.70){\special{em:lineto}}
\put(44.06,21.70){\special{em:moveto}}
\put(44.15,21.70){\special{em:lineto}}
\put(44.40,21.70){\special{em:moveto}}
\put(44.48,21.70){\special{em:lineto}}
\put(44.73,21.70){\special{em:moveto}}
\put(44.81,21.70){\special{em:lineto}}
\put(45.06,21.70){\special{em:moveto}}
\put(45.15,21.70){\special{em:lineto}}
\put(45.40,21.70){\special{em:moveto}}
\put(45.48,21.70){\special{em:lineto}}
\put(45.73,21.70){\special{em:moveto}}
\put(45.81,21.70){\special{em:lineto}}
\put(46.06,21.70){\special{em:moveto}}
\put(46.15,21.70){\special{em:lineto}}
\put(46.40,21.70){\special{em:moveto}}
\put(46.48,21.70){\special{em:lineto}}
\put(46.73,21.70){\special{em:moveto}}
\put(46.81,21.70){\special{em:lineto}}
\put(47.06,21.70){\special{em:moveto}}
\put(47.15,21.70){\special{em:lineto}}
\put(47.40,21.70){\special{em:moveto}}
\put(47.48,21.70){\special{em:lineto}}
\put(47.73,21.70){\special{em:moveto}}
\put(47.81,21.70){\special{em:lineto}}
\put(48.06,21.70){\special{em:moveto}}
\put(48.15,21.70){\special{em:lineto}}
\put(48.40,21.70){\special{em:moveto}}
\put(48.48,21.70){\special{em:lineto}}
\put(48.73,21.70){\special{em:moveto}}
\put(48.81,21.70){\special{em:lineto}}
\put(49.06,21.70){\special{em:moveto}}
\put(49.15,21.70){\special{em:lineto}}
\put(49.40,21.70){\special{em:moveto}}
\put(49.48,21.70){\special{em:lineto}}
\put(49.73,21.70){\special{em:moveto}}
\put(49.81,21.70){\special{em:lineto}}
\put(50.06,21.70){\special{em:moveto}}
\put(50.15,21.70){\special{em:lineto}}
\put(50.40,21.70){\special{em:moveto}}
\put(50.48,21.70){\special{em:lineto}}
\put(50.73,21.70){\special{em:moveto}}
\put(50.81,21.70){\special{em:lineto}}
\put(51.06,21.70){\special{em:moveto}}
\put(51.15,21.70){\special{em:lineto}}
\put(51.40,21.70){\special{em:moveto}}
\put(51.48,21.70){\special{em:lineto}}
\put(51.73,21.70){\special{em:moveto}}
\put(51.81,21.70){\special{em:lineto}}
\put(52.06,21.70){\special{em:moveto}}
\put(52.15,21.70){\special{em:lineto}}
\put(52.40,21.70){\special{em:moveto}}
\put(52.48,21.70){\special{em:lineto}}
\put(52.73,21.70){\special{em:moveto}}
\put(52.81,21.70){\special{em:lineto}}
\put(53.06,21.70){\special{em:moveto}}
\put(53.15,21.70){\special{em:lineto}}
\put(53.40,21.70){\special{em:moveto}}
\put(53.46,21.70){\special{em:lineto}}
\put(23.43,26.66){\special{em:moveto}}
\put(23.51,26.66){\special{em:lineto}}
\put(23.76,26.66){\special{em:moveto}}
\put(23.85,26.66){\special{em:lineto}}
\put(24.10,26.66){\special{em:moveto}}
\put(24.18,26.66){\special{em:lineto}}
\put(24.43,26.66){\special{em:moveto}}
\put(24.51,26.66){\special{em:lineto}}
\put(24.76,26.66){\special{em:moveto}}
\put(24.85,26.66){\special{em:lineto}}
\put(25.10,26.66){\special{em:moveto}}
\put(25.18,26.66){\special{em:lineto}}
\put(25.43,26.66){\special{em:moveto}}
\put(25.51,26.66){\special{em:lineto}}
\put(25.76,26.66){\special{em:moveto}}
\put(25.85,26.66){\special{em:lineto}}
\put(26.10,26.66){\special{em:moveto}}
\put(26.18,26.66){\special{em:lineto}}
\put(26.43,26.66){\special{em:moveto}}
\put(26.51,26.66){\special{em:lineto}}
\put(26.76,26.66){\special{em:moveto}}
\put(26.85,26.66){\special{em:lineto}}
\put(27.10,26.66){\special{em:moveto}}
\put(27.18,26.66){\special{em:lineto}}
\put(27.43,26.66){\special{em:moveto}}
\put(27.51,26.66){\special{em:lineto}}
\put(27.76,26.66){\special{em:moveto}}
\put(27.85,26.66){\special{em:lineto}}
\put(28.10,26.66){\special{em:moveto}}
\put(28.18,26.66){\special{em:lineto}}
\put(28.43,26.66){\special{em:moveto}}
\put(28.51,26.66){\special{em:lineto}}
\put(28.76,26.66){\special{em:moveto}}
\put(28.85,26.66){\special{em:lineto}}
\put(29.10,26.66){\special{em:moveto}}
\put(29.18,26.66){\special{em:lineto}}
\put(29.43,26.66){\special{em:moveto}}
\put(29.51,26.66){\special{em:lineto}}
\put(29.76,26.66){\special{em:moveto}}
\put(29.85,26.66){\special{em:lineto}}
\put(30.10,26.66){\special{em:moveto}}
\put(30.18,26.66){\special{em:lineto}}
\put(30.43,26.66){\special{em:moveto}}
\put(30.51,26.66){\special{em:lineto}}
\put(30.76,26.66){\special{em:moveto}}
\put(30.85,26.66){\special{em:lineto}}
\put(31.10,26.66){\special{em:moveto}}
\put(31.18,26.66){\special{em:lineto}}
\put(31.43,26.66){\special{em:moveto}}
\put(31.51,26.66){\special{em:lineto}}
\put(31.76,26.66){\special{em:moveto}}
\put(31.85,26.66){\special{em:lineto}}
\put(32.10,26.66){\special{em:moveto}}
\put(32.18,26.66){\special{em:lineto}}
\put(32.43,26.66){\special{em:moveto}}
\put(32.51,26.66){\special{em:lineto}}
\put(32.76,26.66){\special{em:moveto}}
\put(32.85,26.66){\special{em:lineto}}
\put(33.10,26.66){\special{em:moveto}}
\put(33.18,26.66){\special{em:lineto}}
\put(33.43,26.66){\special{em:moveto}}
\put(33.51,26.66){\special{em:lineto}}
\put(33.76,26.66){\special{em:moveto}}
\put(33.85,26.66){\special{em:lineto}}
\put(34.10,26.66){\special{em:moveto}}
\put(34.18,26.66){\special{em:lineto}}
\put(34.43,26.66){\special{em:moveto}}
\put(34.51,26.66){\special{em:lineto}}
\put(34.76,26.66){\special{em:moveto}}
\put(34.85,26.66){\special{em:lineto}}
\put(35.10,26.66){\special{em:moveto}}
\put(35.18,26.66){\special{em:lineto}}
\put(35.43,26.66){\special{em:moveto}}
\put(35.51,26.66){\special{em:lineto}}
\put(35.76,26.66){\special{em:moveto}}
\put(35.85,26.66){\special{em:lineto}}
\put(36.10,26.66){\special{em:moveto}}
\put(36.18,26.66){\special{em:lineto}}
\put(36.43,26.66){\special{em:moveto}}
\put(36.51,26.66){\special{em:lineto}}
\put(36.76,26.66){\special{em:moveto}}
\put(36.85,26.66){\special{em:lineto}}
\put(37.10,26.66){\special{em:moveto}}
\put(37.18,26.66){\special{em:lineto}}
\put(37.43,26.66){\special{em:moveto}}
\put(37.51,26.66){\special{em:lineto}}
\put(37.76,26.66){\special{em:moveto}}
\put(37.85,26.66){\special{em:lineto}}
\put(38.10,26.66){\special{em:moveto}}
\put(38.18,26.66){\special{em:lineto}}
\put(38.43,26.66){\special{em:moveto}}
\put(38.51,26.66){\special{em:lineto}}
\put(38.76,26.66){\special{em:moveto}}
\put(38.85,26.66){\special{em:lineto}}
\put(39.10,26.66){\special{em:moveto}}
\put(39.18,26.66){\special{em:lineto}}
\put(39.43,26.66){\special{em:moveto}}
\put(39.51,26.66){\special{em:lineto}}
\put(39.76,26.66){\special{em:moveto}}
\put(39.85,26.66){\special{em:lineto}}
\put(40.10,26.66){\special{em:moveto}}
\put(40.18,26.66){\special{em:lineto}}
\put(40.43,26.66){\special{em:moveto}}
\put(40.51,26.66){\special{em:lineto}}
\put(40.76,26.66){\special{em:moveto}}
\put(40.85,26.66){\special{em:lineto}}
\put(41.10,26.66){\special{em:moveto}}
\put(41.18,26.66){\special{em:lineto}}
\put(41.43,26.66){\special{em:moveto}}
\put(41.51,26.66){\special{em:lineto}}
\put(41.76,26.66){\special{em:moveto}}
\put(41.85,26.66){\special{em:lineto}}
\put(42.10,26.66){\special{em:moveto}}
\put(42.18,26.66){\special{em:lineto}}
\put(42.43,26.66){\special{em:moveto}}
\put(42.51,26.66){\special{em:lineto}}
\put(42.76,26.66){\special{em:moveto}}
\put(42.85,26.66){\special{em:lineto}}
\put(43.10,26.66){\special{em:moveto}}
\put(43.18,26.66){\special{em:lineto}}
\put(43.43,26.66){\special{em:moveto}}
\put(43.51,26.66){\special{em:lineto}}
\put(43.76,26.66){\special{em:moveto}}
\put(43.85,26.66){\special{em:lineto}}
\put(44.10,26.66){\special{em:moveto}}
\put(44.18,26.66){\special{em:lineto}}
\put(44.43,26.66){\special{em:moveto}}
\put(44.51,26.66){\special{em:lineto}}
\put(44.76,26.66){\special{em:moveto}}
\put(44.85,26.66){\special{em:lineto}}
\put(45.10,26.66){\special{em:moveto}}
\put(45.18,26.66){\special{em:lineto}}
\put(45.43,26.66){\special{em:moveto}}
\put(45.51,26.66){\special{em:lineto}}
\put(45.76,26.66){\special{em:moveto}}
\put(45.85,26.66){\special{em:lineto}}
\put(46.10,26.66){\special{em:moveto}}
\put(46.18,26.66){\special{em:lineto}}
\put(46.43,26.66){\special{em:moveto}}
\put(46.51,26.66){\special{em:lineto}}
\put(46.76,26.66){\special{em:moveto}}
\put(46.85,26.66){\special{em:lineto}}
\put(47.10,26.66){\special{em:moveto}}
\put(47.18,26.66){\special{em:lineto}}
\put(47.43,26.66){\special{em:moveto}}
\put(47.51,26.66){\special{em:lineto}}
\put(47.76,26.66){\special{em:moveto}}
\put(47.85,26.66){\special{em:lineto}}
\put(48.10,26.66){\special{em:moveto}}
\put(48.18,26.66){\special{em:lineto}}
\put(48.43,26.66){\special{em:moveto}}
\put(48.51,26.66){\special{em:lineto}}
\put(48.76,26.66){\special{em:moveto}}
\put(48.85,26.66){\special{em:lineto}}
\put(49.10,26.66){\special{em:moveto}}
\put(49.18,26.66){\special{em:lineto}}
\put(49.43,26.66){\special{em:moveto}}
\put(49.51,26.66){\special{em:lineto}}
\put(49.76,26.66){\special{em:moveto}}
\put(49.85,26.66){\special{em:lineto}}
\put(50.10,26.66){\special{em:moveto}}
\put(50.18,26.66){\special{em:lineto}}
\put(50.43,26.66){\special{em:moveto}}
\put(50.51,26.66){\special{em:lineto}}
\put(50.76,26.66){\special{em:moveto}}
\put(50.85,26.66){\special{em:lineto}}
\put(51.10,26.66){\special{em:moveto}}
\put(51.18,26.66){\special{em:lineto}}
\put(51.43,26.66){\special{em:moveto}}
\put(51.51,26.66){\special{em:lineto}}
\put(51.76,26.66){\special{em:moveto}}
\put(51.85,26.66){\special{em:lineto}}
\put(52.10,26.66){\special{em:moveto}}
\put(52.18,26.66){\special{em:lineto}}
\put(52.43,26.66){\special{em:moveto}}
\put(52.51,26.66){\special{em:lineto}}
\put(52.76,26.68){\special{em:moveto}}
\put(52.85,26.68){\special{em:lineto}}
\put(53.10,26.70){\special{em:moveto}}
\put(53.18,26.70){\special{em:lineto}}
\put(26.76,28.36){\special{em:moveto}}
\put(28.43,28.36){\special{em:lineto}}
\put(28.43,30.06){\special{em:moveto}}
\put(28.43,28.36){\special{em:lineto}}
\put(21.86,35.03){\special{em:moveto}}
\put(28.43,28.36){\special{em:lineto}}
\put(56.80,16.70){\special{em:moveto}}
\put(54.93,16.70){\special{em:lineto}}
\put(55.06,18.36){\special{em:moveto}}
\put(55.06,16.66){\special{em:lineto}}
\put(61.73,23.20){\special{em:moveto}}
\put(55.06,16.66){\special{em:lineto}}
\end{picture}
\end{equation}
\end{figure}
  Such a decomposition is not unique: there are exactly $2^{N_i}$ 
  Young diagram's pairs, which lead to the same resulting picture.

  Let us collect all above information in

  \Th\label{nrdata}
  Every connected component $Y\subset\M(2,-1,k)^T$ has a form
  $$\underbrace{\PP_1\times\PP_1\times\cdots\times\PP_1}_{N}\;, 
  $$ 
  and there exists a surjective map from the set of all combinatorial
  data consisting of:

  \smallskip\begingroup\leftskip=1cm

  \noindent\llap{--}
  an ordered triple of non-negative integers $d_0,d_1,d_2$ 
  such that the sum $d=d_0+d_1+d_2$ is bounded by $0\le d\le(k-1)$,

  \smallskip

  \noindent\llap{--}
  an ordered triple of positive integers $a_0,a_1,a_2$, which 
  satisfy three triangle inequalities 
  and the equality $(-a_0+a_1+a_2)^2-4a_1a_2=1-4(c_2(E)-d)$,

  \noindent\llap{--}
  an ordered triple of Yong diagram's pairs $(\l_i,\mu_i)$ such that 
  $|\l_i|+|\mu_i|=d_i$ 

  \smallskip\endgroup

  \noindent
  onto the set of all connected components $Y\subset\M^T$. 
  For any $Y$ there are exactly $2^{\dim Y}$ combinatorial data
  collections, which are mapped into $Y$. They all have the same numbers 
  $a_i$, $d_i$ and the same triples of pictures obtained from the Young
  diagrams in the way shown in (\ref{picsum}). $\PP_1$-multipliers of $Y$
  naturally correspond to the bounded 1-dimensional zones in three diagrams
  (\ref{bfpic}).
\EF

\EP

\No{Some calculations via Bott formula}
\SubNo{Bott formula in our framework}
  Suppose that a torus $T$ acts on a smooth algebraic variety $X$, $\E$ 
  is a toric bundle on $X$, and $\rk(\E)=\dim(X)$. Then the top Chern 
  class of $\E$ can be evaluated only looking on the restriction 
  $\E|_{X^T}$ of $\E$ onto the fixed points loci $X^T$. 

  Namely, let $X^T = \coprod Y$ be the decomposition of $X^T$ into 
  connected components and $\MAP\g,\CC^*,T$ be a general one-parametric 
  subgroup (such that $X^\g=X^T$). For any connected component $Y$ consider 
  the decompositions
  $$\E \mid_{Y} =\bigoplus\limits_{\chi\in\L(T)} \E^{\chi}
    \qquad{\rm and}\qquad
    \nb XY = \bigoplus\limits_{\chi\in\L(T)}\nb XY^{\chi}
  $$ 
  with respect to the torus action. Let
  $$c(\E^{\chi}) = \prod\limits_i(1 + e_{\chi}^i)
    \qquad{\rm and}\qquad
    c(\nb XY^{\chi})= \prod\limits_j (1 + n_{\chi}^j)
  $$ 
  be the formal factorizations of total Chern classes of eigenbundles
   $\E^\chi$, $\nb XY^\chi$.  It follows from the general Bott 
  residue formula
  (see \cite{Bt}) that the top Chern class $\ct(\E)$ can be evaluated
  by the following rule
  $$\ct(\E) = \sum_{Y} \int\limits_{Y} 
    \frac{\prod\limits_{i,\chi}(\sp\g,\chi + e_{\chi}^i)}
         {\prod\limits_{j,\chi}(\sp\g,\chi + n_{\chi}^j)}\quad,
  $$
  where the integrand expression is considered as an element of the Chow 
  ring $A^*(Y)$.

  We are going to apply this formula to calculate $\ct(\G^{\oplus4})$, where
  $\G$ is the  universal bundle over $\M=\M(2,-1,k)$. As we have seen in 
  the previous paragraph, each connected fixed locus $Y\in\M^T$ 
  is represented by a family of all torsion free toric sheaves $E$ 
  on $\PP^2$, which have a given combinatorial type and are equipped 
  with some special toric structures induced by the canonical toric structure
  on $\G$ (see \ref{tsf} above). So, we have to calculate the character's 
  decomposition of $H^1(E,\PP_2)$ and describe the corresponding 
  eigensubbundles of $\G|_Y$ in terms of $A^*(Y)$. Then we need
  to make the same for $\nb{\M(2,-1,k)}Y$, i.~e. we have
  to calculate the character's decomposition of $\ext^1(E,E)$,
  factorize it by the zero character component, and
  describe the corresponding eigensubbundles of $\Ext^1_\pi(\G,\G)|_Y$ 
  in terms of $A^*(Y)$.

\SubNo{Character's decomposition for $\G|_Y$}
  Let $E$ runs through the connected component $Y\subset\M^T$ given in
  combinatorial terms of \ref{nrdata}. Since 
  $H^1(E)=H^1(E^{**})\oplus H^0(C_E)$, 
  the diagram of characters, which  are present in $H^1(E)$, is the sum of
  seven pices: the diagram for $H^1(E^{**})$ and six Young diagrams 
  coming from $H^0(C_i)$, $i=0,1,2$ as it was explained in \ref{hdecomp}
  and \ref{ydpairs}.

  The first pice is calculated directly by looking on the $\chi$-component 
  of the \v{C}hech complex associated with the standard affine 
  covering $\PP_2=\bigcup\limits_{i=0}^2U_i$ (see \cite{K1}).
  It is easy to see that in terms of filtrations $V^{jk}(*)$ described in
  \ref{filttriple}, the character $\chi$ appears in $H^1(E^{**})$ iff
  all three spaces $V^{jk}(\chi)$ ($(jk)=(01),(02),(12)$)
  are 1-dimensional. The set of such characters is represented in $\L(T)$
  by intersection of three strips like
$$
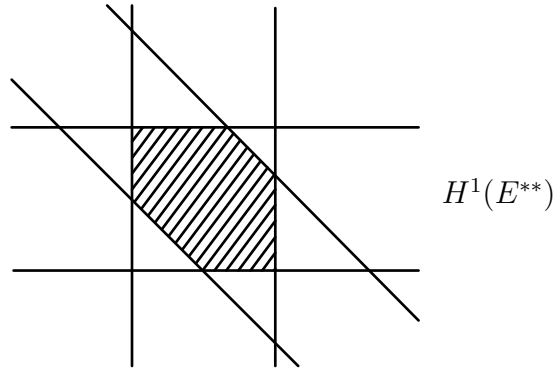
\setlength{\unitlength}{0.075in}
\begin{picture}(37.96,28.36)
\special{em:linewidth 0.014in}
\put(0.00,26.20){\special{em:moveto}}
\put(20.00,6.20){\special{em:lineto}}
\put(6.66,31.36){\special{em:moveto}}
\put(28.40,9.36){\special{em:lineto}}
\put(0.13,12.86){\special{em:moveto}}
\put(28.40,12.86){\special{em:lineto}}
\put(0.00,22.86){\special{em:moveto}}
\put(28.40,22.86){\special{em:lineto}}
\put(18.40,31.20){\special{em:moveto}}
\put(18.40,6.20){\special{em:lineto}}
\put(8.40,31.36){\special{em:moveto}}
\put(8.40,6.20){\special{em:lineto}}
\put(8.40,21.86){\special{em:moveto}}
\put(9.15,22.86){\special{em:lineto}}
\put(8.41,20.75){\special{em:moveto}}
\put(9.98,22.86){\special{em:lineto}}
\put(8.41,19.65){\special{em:moveto}}
\put(10.81,22.86){\special{em:lineto}}
\put(8.40,18.53){\special{em:moveto}}
\put(11.65,22.86){\special{em:lineto}}
\put(8.58,17.66){\special{em:moveto}}
\put(12.48,22.86){\special{em:lineto}}
\put(9.05,17.18){\special{em:moveto}}
\put(13.31,22.86){\special{em:lineto}}
\put(9.53,16.70){\special{em:moveto}}
\put(14.15,22.86){\special{em:lineto}}
\put(10.01,16.23){\special{em:moveto}}
\put(14.98,22.86){\special{em:lineto}}
\put(10.48,15.75){\special{em:moveto}}
\put(15.48,22.43){\special{em:lineto}}
\put(10.95,15.26){\special{em:moveto}}
\put(15.96,21.95){\special{em:lineto}}
\put(11.43,14.80){\special{em:moveto}}
\put(16.43,21.46){\special{em:lineto}}
\put(11.91,14.33){\special{em:moveto}}
\put(16.91,21.00){\special{em:lineto}}
\put(12.38,13.85){\special{em:moveto}}
\put(17.38,20.51){\special{em:lineto}}
\put(12.86,13.36){\special{em:moveto}}
\put(17.85,20.03){\special{em:lineto}}
\put(13.35,12.90){\special{em:moveto}}
\put(18.33,19.55){\special{em:lineto}}
\put(14.15,12.86){\special{em:moveto}}
\put(18.40,18.53){\special{em:lineto}}
\put(14.98,12.86){\special{em:moveto}}
\put(18.40,17.43){\special{em:lineto}}
\put(15.81,12.86){\special{em:moveto}}
\put(18.40,16.33){\special{em:lineto}}
\put(16.65,12.86){\special{em:moveto}}
\put(18.40,15.20){\special{em:lineto}}
\put(17.48,12.86){\special{em:moveto}}
\put(18.40,14.10){\special{em:lineto}}
\put(18.31,12.86){\special{em:moveto}}
\put(18.40,13.00){\special{em:lineto}}
\put(37.95,17.00){{\setbox0=\hbox{$H^1(E^{**})$}\kern-\wd0\raise\dp0\box0}}
\end{picture}
$$

\vspace{-1.5cm}

\noindent
  and all the characters have the multiplicity 1. The weights of the slanted, 
  horizontal, and vertical strips are equal to $a_0$, $a_1$, $a_2$ 
  correspondingly. The same picture for $E^{**}(-1)$ is obtained by the one 
  step shift of the slanting
  strip in the right-upper direction like
 
 $$
  \setlength{\unitlength}{0.075in}
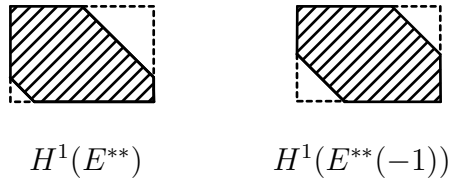
\begin{picture}(30.71,17.66)
\put(19,9){{\setbox0=\hbox{$H^1(E^{**}(-1))$}\lower\ht0\box0}}
\put(2,9){{\setbox0=\hbox{$H^1(E^{**})$}\lower\ht0\box0}}
\special{em:linewidth 0.014in}
\put(0.73,18.23){\special{em:moveto}}
\put(1.11,18.61){\special{em:lineto}}
\put(0.73,17.30){\special{em:moveto}}
\put(2.05,18.61){\special{em:lineto}}
\put(0.73,16.36){\special{em:moveto}}
\put(2.98,18.61){\special{em:lineto}}
\put(0.73,15.43){\special{em:moveto}}
\put(3.91,18.61){\special{em:lineto}}
\put(0.73,14.51){\special{em:moveto}}
\put(4.85,18.61){\special{em:lineto}}
\put(0.78,13.63){\special{em:moveto}}
\put(5.71,18.56){\special{em:lineto}}
\put(1.25,13.16){\special{em:moveto}}
\put(6.18,18.10){\special{em:lineto}}
\put(1.71,12.70){\special{em:moveto}}
\put(6.65,17.63){\special{em:lineto}}
\put(2.18,12.23){\special{em:moveto}}
\put(7.11,17.16){\special{em:lineto}}
\put(2.91,12.03){\special{em:moveto}}
\put(7.58,16.70){\special{em:lineto}}
\put(3.85,12.03){\special{em:moveto}}
\put(8.05,16.23){\special{em:lineto}}
\put(4.78,12.03){\special{em:moveto}}
\put(8.51,15.76){\special{em:lineto}}
\put(5.70,12.03){\special{em:moveto}}
\put(8.98,15.30){\special{em:lineto}}
\put(6.63,12.03){\special{em:moveto}}
\put(9.45,14.85){\special{em:lineto}}
\put(7.56,12.03){\special{em:moveto}}
\put(9.91,14.38){\special{em:lineto}}
\put(8.50,12.03){\special{em:moveto}}
\put(10.38,13.91){\special{em:lineto}}
\put(9.43,12.03){\special{em:moveto}}
\put(10.61,13.21){\special{em:lineto}}
\put(10.36,12.03){\special{em:moveto}}
\put(10.63,12.30){\special{em:lineto}}
\put(0.68,18.66){\special{em:moveto}}
\put(5.68,18.66){\special{em:lineto}}
\put(10.68,13.70){\special{em:lineto}}
\put(10.68,12.83){\special{em:lineto}}
\put(10.71,12.00){\special{em:lineto}}
\put(2.35,12.00){\special{em:lineto}}
\put(0.68,13.66){\special{em:lineto}}
\put(0.68,18.66){\special{em:lineto}}
\put(0.68,13.66){\special{em:moveto}}
\put(0.68,13.33){\special{em:lineto}}
\put(0.68,13.00){\special{em:moveto}}
\put(0.68,12.66){\special{em:lineto}}
\put(0.68,12.33){\special{em:moveto}}
\put(0.68,12.00){\special{em:lineto}}
\put(1.01,12.00){\special{em:moveto}}
\put(1.35,12.00){\special{em:lineto}}
\put(1.68,12.00){\special{em:moveto}}
\put(2.01,12.00){\special{em:lineto}}
\put(2.35,12.00){\special{em:moveto}}
\put(2.38,12.00){\special{em:lineto}}
\put(0.65,13.70){\special{em:moveto}}
\put(2.31,12.03){\special{em:lineto}}
\put(5.68,18.66){\special{em:moveto}}
\put(6.01,18.66){\special{em:lineto}}
\put(6.35,18.66){\special{em:moveto}}
\put(6.68,18.66){\special{em:lineto}}
\put(7.01,18.66){\special{em:moveto}}
\put(7.35,18.66){\special{em:lineto}}
\put(7.68,18.66){\special{em:moveto}}
\put(8.01,18.66){\special{em:lineto}}
\put(8.35,18.66){\special{em:moveto}}
\put(8.68,18.66){\special{em:lineto}}
\put(9.01,18.66){\special{em:moveto}}
\put(9.35,18.66){\special{em:lineto}}
\put(9.68,18.66){\special{em:moveto}}
\put(10.01,18.66){\special{em:lineto}}
\put(10.35,18.66){\special{em:moveto}}
\put(10.68,18.66){\special{em:lineto}}
\put(10.68,18.33){\special{em:moveto}}
\put(10.68,18.00){\special{em:lineto}}
\put(10.68,17.66){\special{em:moveto}}
\put(10.68,17.33){\special{em:lineto}}
\put(10.68,17.00){\special{em:moveto}}
\put(10.68,16.66){\special{em:lineto}}
\put(10.68,16.33){\special{em:moveto}}
\put(10.68,16.00){\special{em:lineto}}
\put(10.68,15.66){\special{em:moveto}}
\put(10.68,15.33){\special{em:lineto}}
\put(10.68,15.00){\special{em:moveto}}
\put(10.68,14.66){\special{em:lineto}}
\put(10.68,14.33){\special{em:moveto}}
\put(10.68,14.00){\special{em:lineto}}
\put(5.68,18.66){\special{em:moveto}}
\put(10.68,13.70){\special{em:lineto}}
\put(2.41,12.00){\special{em:moveto}}
\put(10.68,12.00){\special{em:lineto}}
\put(10.68,13.66){\special{em:lineto}}
\put(20.75,18.33){\special{em:moveto}}
\put(21.06,18.63){\special{em:lineto}}
\put(20.75,17.38){\special{em:moveto}}
\put(22.00,18.63){\special{em:lineto}}
\put(20.75,16.45){\special{em:moveto}}
\put(22.93,18.63){\special{em:lineto}}
\put(20.73,15.50){\special{em:moveto}}
\put(23.86,18.63){\special{em:lineto}}
\put(21.10,14.95){\special{em:moveto}}
\put(24.80,18.63){\special{em:lineto}}
\put(21.56,14.48){\special{em:moveto}}
\put(25.73,18.63){\special{em:lineto}}
\put(22.03,14.01){\special{em:moveto}}
\put(26.65,18.63){\special{em:lineto}}
\put(22.50,13.55){\special{em:moveto}}
\put(27.43,18.48){\special{em:lineto}}
\put(22.96,13.08){\special{em:moveto}}
\put(27.90,18.01){\special{em:lineto}}
\put(23.43,12.61){\special{em:moveto}}
\put(28.36,17.55){\special{em:lineto}}
\put(23.90,12.15){\special{em:moveto}}
\put(28.83,17.08){\special{em:lineto}}
\put(24.71,12.03){\special{em:moveto}}
\put(29.30,16.61){\special{em:lineto}}
\put(25.65,12.03){\special{em:moveto}}
\put(29.76,16.15){\special{em:lineto}}
\put(26.56,12.03){\special{em:moveto}}
\put(30.23,15.68){\special{em:lineto}}
\put(27.50,12.03){\special{em:moveto}}
\put(30.65,15.16){\special{em:lineto}}
\put(28.43,12.03){\special{em:moveto}}
\put(30.65,14.25){\special{em:lineto}}
\put(29.36,12.03){\special{em:moveto}}
\put(30.66,13.31){\special{em:lineto}}
\put(30.30,12.03){\special{em:moveto}}
\put(30.66,12.40){\special{em:lineto}}
\put(20.71,18.66){\special{em:moveto}}
\put(27.31,18.66){\special{em:lineto}}
\put(30.68,15.30){\special{em:lineto}}
\put(30.71,12.00){\special{em:lineto}}
\put(24.01,12.00){\special{em:lineto}}
\put(20.68,15.33){\special{em:lineto}}
\put(20.71,18.66){\special{em:lineto}}
\put(20.68,15.33){\special{em:moveto}}
\put(20.68,15.00){\special{em:lineto}}
\put(20.68,14.66){\special{em:moveto}}
\put(20.68,14.33){\special{em:lineto}}
\put(20.68,14.00){\special{em:moveto}}
\put(20.68,13.66){\special{em:lineto}}
\put(20.68,13.33){\special{em:moveto}}
\put(20.68,13.00){\special{em:lineto}}
\put(20.68,12.66){\special{em:moveto}}
\put(20.68,12.33){\special{em:lineto}}
\put(20.71,12.03){\special{em:moveto}}
\put(21.05,12.03){\special{em:lineto}}
\put(21.38,12.03){\special{em:moveto}}
\put(21.71,12.03){\special{em:lineto}}
\put(22.05,12.03){\special{em:moveto}}
\put(22.38,12.03){\special{em:lineto}}
\put(22.71,12.03){\special{em:moveto}}
\put(23.05,12.03){\special{em:lineto}}
\put(23.38,12.03){\special{em:moveto}}
\put(23.71,12.03){\special{em:lineto}}
\put(24.05,12.03){\special{em:moveto}}
\put(24.15,12.03){\special{em:lineto}}
\put(20.68,15.33){\special{em:moveto}}
\put(24.01,12.00){\special{em:lineto}}
\put(27.35,18.66){\special{em:moveto}}
\put(27.68,18.66){\special{em:lineto}}
\put(28.01,18.66){\special{em:moveto}}
\put(28.35,18.66){\special{em:lineto}}
\put(28.68,18.66){\special{em:moveto}}
\put(29.01,18.66){\special{em:lineto}}
\put(29.35,18.66){\special{em:moveto}}
\put(29.68,18.66){\special{em:lineto}}
\put(30.01,18.66){\special{em:moveto}}
\put(30.35,18.66){\special{em:lineto}}
\put(30.68,18.66){\special{em:moveto}}
\put(30.68,18.33){\special{em:lineto}}
\put(30.68,18.00){\special{em:moveto}}
\put(30.68,17.66){\special{em:lineto}}
\put(30.68,17.33){\special{em:moveto}}
\put(30.68,17.00){\special{em:lineto}}
\put(30.68,16.66){\special{em:moveto}}
\put(30.68,16.33){\special{em:lineto}}
\put(30.68,16.00){\special{em:moveto}}
\put(30.68,15.66){\special{em:lineto}}
\end{picture}
$$

\vspace{-1.5cm}

  Three Young diagram's pairs are placed into three pairs of outward angles
  of the above hexagons in such a way that in more symmetric trigonal 
  coordinates on $\L(T)$ we get the picture shown on (\ref{pichex}).   
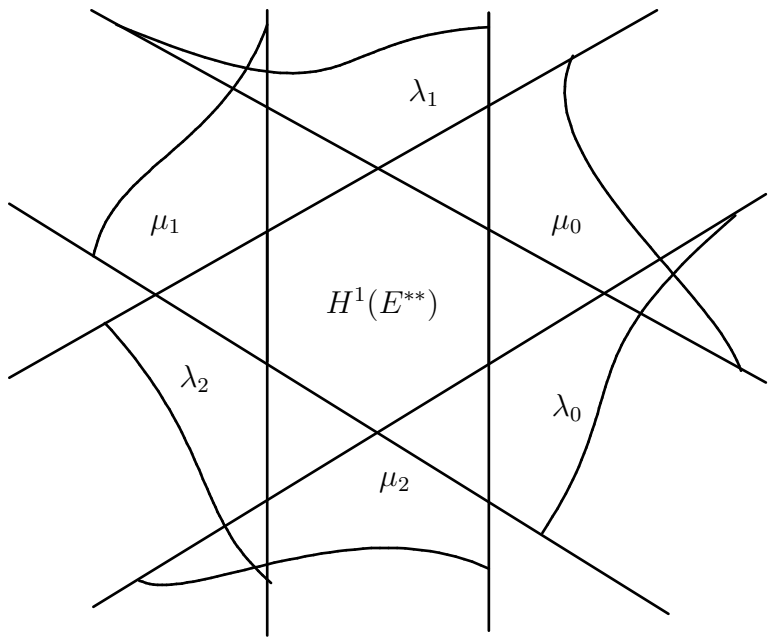
  \begin{equation}\label{pichex} 
\setlength{\unitlength}{0.075in}
\begin{picture}(52.80,43.66)
\special{em:linewidth 0.014in}
\put(30,22){{\setbox0=\hbox{$H^1(E^{**})$}\kern-\wd0\raise\dp0\box0}}
\put(28,10){{\setbox0=\hbox{$\mu_2$}\kern-\wd0\raise\dp0\box0}}
\put(14,17){{\setbox0=\hbox{$\lambda_2$}\kern-\wd0\raise\dp0\box0}}
\put(12,28){{\setbox0=\hbox{$\mu_1$}\kern-\wd0\raise\dp0\box0}}
\put(30,37){{\setbox0=\hbox{$\lambda_1$}\kern-\wd0\raise\dp0\box0}}
\put(40,28){{\setbox0=\hbox{$\mu_0$}\kern-\wd0\raise\dp0\box0}}
\put(40,15){{\setbox0=\hbox{$\lambda_0$}\kern-\wd0\raise\dp0\box0}}
\put(33.46,42.50){\special{em:moveto}}
\put(32.75,42.45){\special{em:lineto}}
\put(32.06,42.38){\special{em:lineto}}
\put(31.43,42.31){\special{em:lineto}}
\put(30.83,42.23){\special{em:lineto}}
\put(30.25,42.13){\special{em:lineto}}
\put(29.70,42.03){\special{em:lineto}}
\put(29.18,41.93){\special{em:lineto}}
\put(28.68,41.81){\special{em:lineto}}
\put(28.21,41.70){\special{em:lineto}}
\put(27.76,41.56){\special{em:lineto}}
\put(27.33,41.43){\special{em:lineto}}
\put(26.93,41.30){\special{em:lineto}}
\put(26.53,41.16){\special{em:lineto}}
\put(26.15,41.01){\special{em:lineto}}
\put(25.78,40.88){\special{em:lineto}}
\put(25.43,40.73){\special{em:lineto}}
\put(25.08,40.60){\special{em:lineto}}
\put(24.73,40.46){\special{em:lineto}}
\put(24.40,40.33){\special{em:lineto}}
\put(24.06,40.20){\special{em:lineto}}
\put(23.75,40.06){\special{em:lineto}}
\put(23.41,39.95){\special{em:lineto}}
\put(23.08,39.83){\special{em:lineto}}
\put(22.76,39.73){\special{em:lineto}}
\put(22.43,39.63){\special{em:lineto}}
\put(22.08,39.55){\special{em:lineto}}
\put(21.73,39.46){\special{em:lineto}}
\put(21.38,39.40){\special{em:lineto}}
\put(21.00,39.35){\special{em:lineto}}
\put(20.61,39.31){\special{em:lineto}}
\put(20.21,39.28){\special{em:lineto}}
\put(19.80,39.26){\special{em:lineto}}
\put(19.36,39.28){\special{em:lineto}}
\put(18.91,39.30){\special{em:lineto}}
\put(18.45,39.33){\special{em:lineto}}
\put(17.95,39.40){\special{em:lineto}}
\put(17.41,39.46){\special{em:lineto}}
\put(16.86,39.56){\special{em:lineto}}
\put(16.28,39.68){\special{em:lineto}}
\put(15.66,39.83){\special{em:lineto}}
\put(15.03,39.98){\special{em:lineto}}
\put(14.35,40.18){\special{em:lineto}}
\put(13.63,40.38){\special{em:lineto}}
\put(12.88,40.63){\special{em:lineto}}
\put(12.08,40.90){\special{em:lineto}}
\put(11.25,41.18){\special{em:lineto}}
\put(10.38,41.51){\special{em:lineto}}
\put(9.45,41.86){\special{em:lineto}}
\put(8.48,42.25){\special{em:lineto}}
\put(7.46,42.66){\special{em:lineto}}
\put(5.86,26.50){\special{em:moveto}}
\put(5.95,26.86){\special{em:lineto}}
\put(6.05,27.23){\special{em:lineto}}
\put(6.16,27.58){\special{em:lineto}}
\put(6.28,27.91){\special{em:lineto}}
\put(6.41,28.25){\special{em:lineto}}
\put(6.56,28.56){\special{em:lineto}}
\put(6.73,28.88){\special{em:lineto}}
\put(6.90,29.18){\special{em:lineto}}
\put(7.10,29.48){\special{em:lineto}}
\put(7.28,29.76){\special{em:lineto}}
\put(7.50,30.05){\special{em:lineto}}
\put(7.71,30.33){\special{em:lineto}}
\put(7.93,30.60){\special{em:lineto}}
\put(8.18,30.86){\special{em:lineto}}
\put(8.41,31.13){\special{em:lineto}}
\put(8.66,31.40){\special{em:lineto}}
\put(8.93,31.65){\special{em:lineto}}
\put(9.20,31.91){\special{em:lineto}}
\put(9.46,32.16){\special{em:lineto}}
\put(9.75,32.41){\special{em:lineto}}
\put(10.01,32.68){\special{em:lineto}}
\put(10.31,32.93){\special{em:lineto}}
\put(10.60,33.18){\special{em:lineto}}
\put(10.90,33.43){\special{em:lineto}}
\put(11.18,33.70){\special{em:lineto}}
\put(11.48,33.96){\special{em:lineto}}
\put(11.80,34.23){\special{em:lineto}}
\put(12.10,34.50){\special{em:lineto}}
\put(12.40,34.76){\special{em:lineto}}
\put(12.70,35.05){\special{em:lineto}}
\put(13.00,35.33){\special{em:lineto}}
\put(13.30,35.61){\special{em:lineto}}
\put(13.60,35.91){\special{em:lineto}}
\put(13.90,36.23){\special{em:lineto}}
\put(14.20,36.53){\special{em:lineto}}
\put(14.50,36.86){\special{em:lineto}}
\put(14.78,37.20){\special{em:lineto}}
\put(15.06,37.53){\special{em:lineto}}
\put(15.35,37.90){\special{em:lineto}}
\put(15.63,38.25){\special{em:lineto}}
\put(15.90,38.63){\special{em:lineto}}
\put(16.16,39.01){\special{em:lineto}}
\put(16.43,39.41){\special{em:lineto}}
\put(16.68,39.83){\special{em:lineto}}
\put(16.91,40.26){\special{em:lineto}}
\put(17.16,40.71){\special{em:lineto}}
\put(17.38,41.18){\special{em:lineto}}
\put(17.60,41.65){\special{em:lineto}}
\put(17.81,42.15){\special{em:lineto}}
\put(18.00,42.66){\special{em:lineto}}
\put(6.66,21.83){\special{em:moveto}}
\put(7.13,21.33){\special{em:lineto}}
\put(7.56,20.81){\special{em:lineto}}
\put(7.98,20.31){\special{em:lineto}}
\put(8.38,19.83){\special{em:lineto}}
\put(8.76,19.35){\special{em:lineto}}
\put(9.11,18.88){\special{em:lineto}}
\put(9.46,18.41){\special{em:lineto}}
\put(9.78,17.96){\special{em:lineto}}
\put(10.08,17.51){\special{em:lineto}}
\put(10.38,17.06){\special{em:lineto}}
\put(10.65,16.63){\special{em:lineto}}
\put(10.91,16.21){\special{em:lineto}}
\put(11.16,15.80){\special{em:lineto}}
\put(11.40,15.38){\special{em:lineto}}
\put(11.61,14.98){\special{em:lineto}}
\put(11.83,14.58){\special{em:lineto}}
\put(12.03,14.18){\special{em:lineto}}
\put(12.21,13.80){\special{em:lineto}}
\put(12.40,13.41){\special{em:lineto}}
\put(12.58,13.05){\special{em:lineto}}
\put(12.75,12.68){\special{em:lineto}}
\put(12.91,12.31){\special{em:lineto}}
\put(13.06,11.95){\special{em:lineto}}
\put(13.23,11.60){\special{em:lineto}}
\put(13.36,11.26){\special{em:lineto}}
\put(13.51,10.91){\special{em:lineto}}
\put(13.66,10.58){\special{em:lineto}}
\put(13.81,10.25){\special{em:lineto}}
\put(13.95,9.91){\special{em:lineto}}
\put(14.10,9.58){\special{em:lineto}}
\put(14.23,9.26){\special{em:lineto}}
\put(14.38,8.95){\special{em:lineto}}
\put(14.53,8.63){\special{em:lineto}}
\put(14.68,8.33){\special{em:lineto}}
\put(14.85,8.01){\special{em:lineto}}
\put(15.00,7.71){\special{em:lineto}}
\put(15.18,7.41){\special{em:lineto}}
\put(15.35,7.11){\special{em:lineto}}
\put(15.53,6.81){\special{em:lineto}}
\put(15.71,6.51){\special{em:lineto}}
\put(15.91,6.23){\special{em:lineto}}
\put(16.13,5.95){\special{em:lineto}}
\put(16.35,5.65){\special{em:lineto}}
\put(16.58,5.36){\special{em:lineto}}
\put(16.83,5.08){\special{em:lineto}}
\put(17.10,4.80){\special{em:lineto}}
\put(17.36,4.51){\special{em:lineto}}
\put(17.65,4.23){\special{em:lineto}}
\put(17.95,3.95){\special{em:lineto}}
\put(18.26,3.66){\special{em:lineto}}
\put(33.46,4.66){\special{em:moveto}}
\put(32.88,4.91){\special{em:lineto}}
\put(32.31,5.15){\special{em:lineto}}
\put(31.73,5.33){\special{em:lineto}}
\put(31.15,5.50){\special{em:lineto}}
\put(30.56,5.65){\special{em:lineto}}
\put(29.98,5.78){\special{em:lineto}}
\put(29.40,5.88){\special{em:lineto}}
\put(28.80,5.96){\special{em:lineto}}
\put(28.21,6.03){\special{em:lineto}}
\put(27.61,6.08){\special{em:lineto}}
\put(27.05,6.10){\special{em:lineto}}
\put(26.45,6.11){\special{em:lineto}}
\put(25.86,6.11){\special{em:lineto}}
\put(25.26,6.10){\special{em:lineto}}
\put(24.68,6.06){\special{em:lineto}}
\put(24.10,6.01){\special{em:lineto}}
\put(23.51,5.95){\special{em:lineto}}
\put(22.93,5.88){\special{em:lineto}}
\put(22.36,5.81){\special{em:lineto}}
\put(21.80,5.73){\special{em:lineto}}
\put(21.23,5.63){\special{em:lineto}}
\put(20.66,5.53){\special{em:lineto}}
\put(20.11,5.41){\special{em:lineto}}
\put(19.56,5.31){\special{em:lineto}}
\put(19.01,5.20){\special{em:lineto}}
\put(18.48,5.08){\special{em:lineto}}
\put(17.96,4.95){\special{em:lineto}}
\put(17.45,4.83){\special{em:lineto}}
\put(16.93,4.71){\special{em:lineto}}
\put(16.43,4.58){\special{em:lineto}}
\put(15.95,4.46){\special{em:lineto}}
\put(15.46,4.35){\special{em:lineto}}
\put(15.00,4.25){\special{em:lineto}}
\put(14.53,4.13){\special{em:lineto}}
\put(14.08,4.03){\special{em:lineto}}
\put(13.65,3.93){\special{em:lineto}}
\put(13.23,3.85){\special{em:lineto}}
\put(12.81,3.76){\special{em:lineto}}
\put(12.41,3.70){\special{em:lineto}}
\put(12.03,3.63){\special{em:lineto}}
\put(11.66,3.60){\special{em:lineto}}
\put(11.31,3.55){\special{em:lineto}}
\put(10.96,3.53){\special{em:lineto}}
\put(10.65,3.53){\special{em:lineto}}
\put(10.35,3.53){\special{em:lineto}}
\put(10.05,3.56){\special{em:lineto}}
\put(9.78,3.60){\special{em:lineto}}
\put(9.53,3.66){\special{em:lineto}}
\put(9.30,3.75){\special{em:lineto}}
\put(9.06,3.83){\special{em:lineto}}
\put(50.66,29.33){\special{em:moveto}}
\put(49.95,28.71){\special{em:lineto}}
\put(49.28,28.11){\special{em:lineto}}
\put(48.65,27.55){\special{em:lineto}}
\put(48.06,27.00){\special{em:lineto}}
\put(47.51,26.46){\special{em:lineto}}
\put(46.98,25.93){\special{em:lineto}}
\put(46.50,25.43){\special{em:lineto}}
\put(46.05,24.93){\special{em:lineto}}
\put(45.61,24.46){\special{em:lineto}}
\put(45.23,24.00){\special{em:lineto}}
\put(44.85,23.55){\special{em:lineto}}
\put(44.51,23.11){\special{em:lineto}}
\put(44.20,22.68){\special{em:lineto}}
\put(43.90,22.26){\special{em:lineto}}
\put(43.63,21.85){\special{em:lineto}}
\put(43.38,21.45){\special{em:lineto}}
\put(43.15,21.05){\special{em:lineto}}
\put(42.93,20.66){\special{em:lineto}}
\put(42.73,20.28){\special{em:lineto}}
\put(42.55,19.90){\special{em:lineto}}
\put(42.38,19.53){\special{em:lineto}}
\put(42.23,19.16){\special{em:lineto}}
\put(42.08,18.80){\special{em:lineto}}
\put(41.95,18.43){\special{em:lineto}}
\put(41.83,18.06){\special{em:lineto}}
\put(41.70,17.71){\special{em:lineto}}
\put(41.58,17.35){\special{em:lineto}}
\put(41.48,16.98){\special{em:lineto}}
\put(41.36,16.61){\special{em:lineto}}
\put(41.26,16.25){\special{em:lineto}}
\put(41.15,15.88){\special{em:lineto}}
\put(41.05,15.50){\special{em:lineto}}
\put(40.93,15.11){\special{em:lineto}}
\put(40.81,14.73){\special{em:lineto}}
\put(40.68,14.35){\special{em:lineto}}
\put(40.55,13.95){\special{em:lineto}}
\put(40.41,13.53){\special{em:lineto}}
\put(40.26,13.11){\special{em:lineto}}
\put(40.11,12.70){\special{em:lineto}}
\put(39.93,12.26){\special{em:lineto}}
\put(39.75,11.81){\special{em:lineto}}
\put(39.55,11.36){\special{em:lineto}}
\put(39.33,10.88){\special{em:lineto}}
\put(39.08,10.40){\special{em:lineto}}
\put(38.83,9.90){\special{em:lineto}}
\put(38.55,9.40){\special{em:lineto}}
\put(38.25,8.86){\special{em:lineto}}
\put(37.91,8.31){\special{em:lineto}}
\put(37.56,7.75){\special{em:lineto}}
\put(37.20,7.16){\special{em:lineto}}
\put(0.00,30.16){\special{em:moveto}}
\put(46.00,1.50){\special{em:lineto}}
\put(18.00,19.16){\special{em:moveto}}
\put(18.00,43.66){\special{em:lineto}}
\put(18.00,28.16){\special{em:moveto}}
\put(18.00,0.00){\special{em:lineto}}
\put(39.33,40.50){\special{em:moveto}}
\put(39.13,40.06){\special{em:lineto}}
\put(38.98,39.63){\special{em:lineto}}
\put(38.88,39.20){\special{em:lineto}}
\put(38.80,38.76){\special{em:lineto}}
\put(38.76,38.33){\special{em:lineto}}
\put(38.75,37.88){\special{em:lineto}}
\put(38.78,37.43){\special{em:lineto}}
\put(38.83,36.98){\special{em:lineto}}
\put(38.91,36.51){\special{em:lineto}}
\put(39.03,36.06){\special{em:lineto}}
\put(39.18,35.60){\special{em:lineto}}
\put(39.33,35.13){\special{em:lineto}}
\put(39.53,34.66){\special{em:lineto}}
\put(39.75,34.20){\special{em:lineto}}
\put(39.98,33.71){\special{em:lineto}}
\put(40.23,33.25){\special{em:lineto}}
\put(40.50,32.78){\special{em:lineto}}
\put(40.80,32.30){\special{em:lineto}}
\put(41.10,31.83){\special{em:lineto}}
\put(41.41,31.35){\special{em:lineto}}
\put(41.75,30.88){\special{em:lineto}}
\put(42.10,30.40){\special{em:lineto}}
\put(42.45,29.93){\special{em:lineto}}
\put(42.81,29.46){\special{em:lineto}}
\put(43.18,29.00){\special{em:lineto}}
\put(43.56,28.51){\special{em:lineto}}
\put(43.95,28.05){\special{em:lineto}}
\put(44.35,27.58){\special{em:lineto}}
\put(44.73,27.13){\special{em:lineto}}
\put(45.11,26.66){\special{em:lineto}}
\put(45.51,26.21){\special{em:lineto}}
\put(45.90,25.75){\special{em:lineto}}
\put(46.28,25.30){\special{em:lineto}}
\put(46.66,24.86){\special{em:lineto}}
\put(47.05,24.41){\special{em:lineto}}
\put(47.41,23.98){\special{em:lineto}}
\put(47.76,23.55){\special{em:lineto}}
\put(48.11,23.13){\special{em:lineto}}
\put(48.46,22.71){\special{em:lineto}}
\put(48.78,22.30){\special{em:lineto}}
\put(49.10,21.88){\special{em:lineto}}
\put(49.40,21.48){\special{em:lineto}}
\put(49.68,21.08){\special{em:lineto}}
\put(49.93,20.70){\special{em:lineto}}
\put(50.18,20.31){\special{em:lineto}}
\put(50.40,19.93){\special{em:lineto}}
\put(50.60,19.56){\special{em:lineto}}
\put(50.78,19.20){\special{em:lineto}}
\put(50.93,18.85){\special{em:lineto}}
\put(51.06,18.50){\special{em:lineto}}
\put(45.20,43.66){\special{em:moveto}}
\put(0.00,18.00){\special{em:lineto}}
\put(52.80,30.83){\special{em:moveto}}
\put(5.86,2.00){\special{em:lineto}}
\put(5.73,43.66){\special{em:moveto}}
\put(52.80,17.66){\special{em:lineto}}
\put(33.46,43.50){\special{em:moveto}}
\put(33.46,0.33){\special{em:lineto}}
\put(33.46,19.00){\special{em:moveto}}
\put(33.46,43.50){\special{em:lineto}}
\end{picture}
  \end{equation} 
  Recall that this picture must be placed into $\L(T)$ in such a way that
  the sum of all characters from $H^1(E)$ equals to the sum of all 
  characters from $H^1(E(-1))$.

  Eigensubbundle corresponding to the character $\chi$ is either trivial or
  the pull-back of $\O_{\PP_1}(1)$. The last case was explained in 
  \ref{hdecomp}.

\SubNo{Character's decomposition for $\nb{\M}Y$}
  We compute the character decomposition of $\T_{\M,E}=Ext^1(E,E)$ using
  the standard spectral sequence associated with the monad (\ref{m}). In 
  this sequence $E_2^{pq}=E_\infty^{pq}$ and $E_1^{pq}$ coincides 
  with the following $T$-equvariant complex of $T$-modules
  \begin{equation}\label{monforext}
   0\TO{}\begin{array}{c}
            \End(U_0)\\ \oplus\\ \End(U_1)\\ \oplus\\ \End(U_2)
          \end{array}
     \TO{}\V\otimes\begin{array}{c}
            \hom(H^1(U_0,U_1)\\ \oplus\\ 
            \hom(H^1(U_1,U_2)
          \end{array} 
     \TO{}\L^2(\V)\otimes\hom(U_0,U_2)
     \TO{}0
  \end{equation}
   where $U_i=H^1(E(i-2))$. This complex has two cohomology spaces: 
  $\CC$ in the left term and $\ext^1(E,E)$ in the middle term. Hence, 
  it gives the character table for $\T_{\M,E}$. Unfortunately, we
  can not represent the answer by a nice picture like in the previous 
  section. But the calculation has the strightforward 
  algorithmization in terms of
  set-theoretical and Minkowski sums of character tables. Since we have
  a full description of all $T$-modules from (\ref{monforext}), we can
  calculate the character table for $\ext^1(E,E)$ by computer.
  The character decomposition of $\nb \M{Y}$ can be extracted from this 
  table immediately by omitting the zero character. Moreover, using 
  the formal decomposition of Chern polynomials for eigensubbundles in 
  $H^1(E(i))$ presented in \ref{hdecomp}, we can extract from 
  (\ref{monforext}) not only the character table but also the formal
  expansions for the denominators in the Bott formula. We can make this 
  last step also only numerically.

  Exact numerical results obtained by this algorithm are the following:

  $$\mbox{\begin{tabular}{|c||c|c|c|c|c|c|c|}
              \hline
              $k$&1&2&3&4&5&6&7\\
              \hline 
              $a_k$&0&0&0&13&729&85026&15650066\\
              \hline 
          \end{tabular} 
         }\;.
  $$

\end{document}